\documentclass{article}

\usepackage{arxiv}

\usepackage[utf8]{inputenc} 
\usepackage[T1]{fontenc}    
\usepackage{hyperref}       
\usepackage{url}            
\usepackage{booktabs}       
\usepackage{amsfonts}       
\usepackage{nicefrac}       
\usepackage{microtype}      
\usepackage{lipsum}
\usepackage[utf8]{inputenc} 
\usepackage[T1]{fontenc}    
\usepackage{hyperref}       
\usepackage{url}           
\usepackage{booktabs}      
\usepackage{amsfonts}     
\usepackage{nicefrac}       
\usepackage{microtype}     
\usepackage{lipsum}
\usepackage{tikz}
\usepackage[noend]{algpseudocode}
\usepackage{algorithmicx, algorithm}
\usepackage{soul, color, xcolor}
\usepackage[misc]{ifsym}
\usepackage{graphicx}
\usepackage{epstopdf}
\usepackage{amsmath}
\usepackage{amsthm}
\usepackage{todonotes}
\usepackage{comment}
\usepackage{cases}
\usepackage{subcaption}

\usepackage{lineno}
\soulregister\ref7
\soulregister\cite7

\theoremstyle{remark}
\newtheorem*{remark}{Remark}

\graphicspath{ {./Figures/} }

\title{A conservative hybrid physics-informed neural network method for {Maxwell--Amp\`{e}re--Nernst--Planck} equations}

\author{
 Cheng Chang \thanks{Department of Mathematics \& Institute of Mathematical Sciences, The Chinese University of Hong Kong. Email: chengchang@link.cuhk.edu.hk. Address: Room G08, Lady Shaw Building, The Chinese University of Hong Kong, Shatin, N.T., Hong Kong Special Administrative Region, People's Republic of China},~~
 Zhouping Xin \thanks{Institute of Mathematical Sciences, The Chinese University of Hong Kong. Email: zpxin@ims.cuhk.edu.hk. Address: Academic Building No. 1, The Chinese University of Hong Kong, Shatin, N.T., Hong Kong Special Administrative Region, People's Republic of China},~~
  Tieyong Zeng \thanks{Department of Mathematics, The Chinese University of Hong Kong. Email: zeng@math.cuhk.edu.hk. Address: Room 225, Lady Shaw Building, The Chinese University of Hong Kong, Shatin, N.T., Hong Kong Special Administrative Region, People's Republic of China}
 }

\begin{document}
\maketitle

\begin{abstract}
{Maxwell--Amp\`{e}re--Nernst--Planck} ({M}ANP) equations {were} recently proposed to model the dynamics of charged particles. In this {study}, we enhance {a} numerical algorithm {of this system} with deep learning tools. The proposed hybrid algorithm provides an automated means to {determine a} proper approximation for the dummy variables, which {can} otherwise only be obtained through massive numerical tests. {In addition}, the original method is validated {for} 2-dimensional problems. However, when the spatial dimension is {one}, the original curl-free relaxation component {is inapplicable}, and the approximation formula for dummy variables{,} which works well in {a} 2-dimensional scenario{,} {fails} to provide a reasonable output {in the 1-dimensional case}. The proposed method {can be} readily generalised to cases {with one spatial dimension}. Experiments show{ }numerical stability and good convergence to the steady{-}state solution obtained from Poisson{--}Boltzmann type equations in {the} 1-dimensional case. The experiments conducted in {the} 2-dimensional case indicate that the proposed method preserves the conservation properties.{ }
\end{abstract}

\keywords{{Maxwell--Amp\`{e}re--Nernst--Planck} equations \and Conservative numerical scheme \and Deep learning \and Physics-informed neural network}

\section{Introduction}
During the past few years, { }deep learning technologies have flourished in {various} applications, such as image processing \cite{Gfm001}, speech recognition \cite{zengSR006}, and natural language processing \cite{FPDL14}. {In addition to} those high-level tasks which {have traditionally thought to} require human's intelligence, deep learning techniques also demonstrate their efficacy in dealing with low-level problems such as partial differential equations (PDEs). Dated back to 1990s, the authors of \cite{PINN006} use the total $l_2$ error corresponding to the governing equations and boundary conditions at collocation points as {the} loss function{. B}y minimising it the neural network is trained to fit the equations. {I}n \cite{UQDL004}{,} the authors propose more sophisticated formulations{ }in which the initial and boundary conditions are enforced by construction{;} thus{,} only the residual for the governing equations need to be included in the training loss. Recently, more algorithms {applying} deep learning to solve PDEs were proposed. In \cite{PINN001}, the authors propose a mesh-free method, and meanwhile provide an efficient approximate algorithm to calculate expressions involving a large number of second-order derivatives, {thereby alleviating} the {c}urse of {d}imensionality \cite{AWS03}. {Elliptic equations can be solved by} optimising the variational problem corresponding to the original equations \cite{EEDL08}. {Building} upon {these} proceeding works, Raissi et al. develop the impactful Physics{--i}nformed {n}eural {n}etwork (PINN){, providing two different paradigms for solving PDEs with PINN, i.e. continuous and discrete time model for both the forward and inverse problems \cite{DLDE001}}. After those cornerstone works, numerous variants on PINN are proposed. To name a few, in \cite{PINN003}, in order to tackle the jump discontinuities which is common in the solutions of {conservation laws} (even when smooth initial conditions are given), the authors propose to approximate the solutions in sub-domains with separate PINNs, then combine them{ }with flux at the sub-domain boundary. In variational PINNs, the authors {utilise} the weak forms of the residual to train{ }neural networks{ }with some set of test functions (e.g.{,} trigonometric functions or polynomials) \cite{PINN007}. For the weak residual of the PDE, {in addition to} simple test functions, \cite{UQDL003} {employ} one neural network to fit the solution, {whereas the other} neural network is used as the test function. As an improvement of VPINN, in \cite{PINN004}, the authors investigate several popular PINN variants{,} then point out that their essential difference is the choice of the test functions in the weak formulations, and propose to use localised test functions to assist the training of the neural network in each sub-domain. 

Despite the efficacy of PINNs {in} solving PDEs, the conventional numerical algorithms are{ }incorporated to remedy some defects in deep learning approaches. In \cite{ANP064}, the authors discover that the deep neural networks fit low-frequency features first, which {differs} from many {conventional} numerical methods. Considering this, a hybrid method is developed, in which the PINN is trained to obtain a coarse-grained solution first, {and} then a {conventional} linear solver is applied to refine the solution \cite{ANP055}. A similar method is proposed in \cite{ANP065}, however, in a reverse way, where a {conventional} multigrid solver is {applied} first to get {the} solution in the coarsest few level{s} and then it comes to {the} deep neural network to rectify the solution to finer levels. Besides using classical numerical schemes and deep learning at different resolution levels, expecting them to perform well in their respective layers, {attempts have been made} to merge them {more interactively}. {Studies} in this direction usually do not {rely entirely} on {a}utomatic {d}ifferentiation (AD) \cite{ANP069}, but {utilise} some {conventional} numerical methods to{ }discretise the domain and approximate differential operators. For example, {the discrete-time model of PINN} has a similar form to classical Runge{--}Kutta methods, except that{ }intermediate time-step values are generated with a neural network {\cite{DLDE001}}. Fang {utilises} well-designed examples to illustrate the faults brought by AD, and then proposes a generalised finite{ }difference-like scheme to replace the AD \cite{ANP033}. { }Finite-{d}ifference (FD) approximations for differential operators are also {adopted} in \cite{ANP047}, {whereas} the cell size is adaptively shrunk near the boundary {such} that the neighbouring grid points are guaranteed to be inside the domain of interest. In \cite{ANP045}, {a} neural network approximation with AD is applied {in} regions where the solution is smooth, {whereas} a {Weighted Essentially Non-oscillatory scheme} is {adopted in} non-smooth regions. In \cite{ANP048}, the authors propose a number of second-order accurate numerical differentiation schemes {in which} the derivative terms are obtained through a {combination} of AD and FD. In \cite{ANP052}, {a} neural network is used as a corrector for partially unknown physics, {and} a splitting scheme in time is {utilised} to tackle {two} sub-problems{;} one with the neural network and the other without it. In \cite{ANP054}, an auto-encoder is trained to obtain{ }latent variable representations, and then the image gradient (a central-difference scheme) operator is used to approximate the second-order derivatives appearing in {the} heat equation. Auto-encoder is also {utilised} in \cite{ANP059} to encode not only the solution field (velocity, pressure in Navier{--}Stokes equations), but also the domain geometry and boundary condition information{. M}eanwhile, the PDE residual is calculated {with} classical finite volume schemes.

One {of the reasons why} the {conventional} numerical methods {are}, in some cases, preferable to deep learning {approaches} is that the solution is guaranteed to fulfil certain physical properties in carefully{ }designed traditional numerical methods. {Taking} the Poisson{--}Nernst{--}Planck (PNP) equations as an example, this system of equations models the dynamics of the ions under {a} background electric field, which{ }is affected by the distribution of those ions. The total amount of each species of {ions} remains constant (provided proper boundary condition, e.g.{,} no flux), while the total energy dissipates as time elapses \cite{NP002}. When a numerical scheme is {carefully designed}, {these} properties can be analytically preserved{ }either conditionally or unconditionally \cite{NP002, NP003, NP008}. {Some studies have attempted to preserve} such kind of physical properties in the context of deep learning. In \cite{ANP067}, a special network structure is constructed to constrain the trained network to be symplectic, thus {better fitting the phase flow of Hamiltonian systems}. In \cite{ANP066}, simple physical constraints {such as} odd and even constraints are guaranteed at specially{ }designed hub layers, and more complicated constraints{, such as energy conservation} are ensured by a separate corrector network. {In addition to the} hard constraints imposed on{ }neural network structures, in \cite{NP015}{,} the authors propose {adding} soft constraints on physical properties (e.g. mass conservation) to the neural network {by introducing a} penalty term into the training loss.

In this {study}, we propose a numerical method which is built upon {conventional} FD discretisation \cite{ANP000_2} and introduce a neural network approximator to solve {Maxwell--}Amp\`{e}re{--}Nernst{--}Planck ({M}ANP) equations \cite{ANP000_1}. As the original scheme{ }, the proposed hybrid scheme also possesses properties {such as mass conservation and positivity preservation}. Furthermore, it does not rely on experimental formulas{ }and is more flexible when the problem settings {are altered}. The {remainder} of this paper is organised as follows. In Section \ref{sec:ANP_base}, we give a brief review on {M}ANP equations and the original numerical scheme. In Section \ref{sec:PINN}, we briefly review the general framework of {the} PINN. We present our proposed hybrid PINN method in Section \ref{sec:hybrid_PINN}. The numerical experiments are in Section \ref{sec:NE}

\section{\texorpdfstring{{Maxwell--}}{Maxwell--}Amp\`{e}re\texorpdfstring{{--}}{--}Nernst\texorpdfstring{{--}}{--}Planck equations and its numerical method}
\label{sec:ANP_base}

{Consider the equations} on spatial domain $\Omega\subset\mathbb{R}^d$ and time interval $[0, T]$\footnote[1]{{For simplicity, we omit those physics constant coefficients and mean field approximation correction terms, which are of little importance for the presentation of the proposed numerical scheme.}}:

\begin{numcases}{}
	\frac{\partial c^l}{\partial t}= -\nabla\cdot\boldsymbol{J}^l\text{, }l=1,2,...,M\label{NP_1}\\ 
	{\frac{\partial \boldsymbol{D}}{\partial t}=-\sum_{l=1}^M q^l\boldsymbol{J}^l+\boldsymbol{\Theta}\label{A_1}}\\
	{\boldsymbol{J}^l=-(\nabla c^l-\frac{q^lc^l}{\epsilon}\boldsymbol{D})\text{, }l=1,2,...,M\label{NP_2_A_2}}\\
	\nabla\cdot\boldsymbol{\Theta}=0\label{Theta}\label{Theta_div_free}\\
	\nabla\times\frac{\boldsymbol{D}}{\epsilon}=0\label{D_curl_free}
\end{numcases}

in which $l$ is the index for the species of ions, $c^l:\Omega\times[0,T]\rightarrow\mathbb{R}^+$ is the concentration of the $l$-th ions,{ }$\boldsymbol{D}:\Omega\times[0,T]\rightarrow\mathbb{R}^d$ is the electric displacement relating to electric potential $\phi$ via $\boldsymbol{D}/\epsilon=-\nabla\phi$, $\epsilon:\Omega\rightarrow\mathbb{R}^+$ is the permittivity, $q^l\in\mathbb{Z}$ is the valence of the $l$-th ions{. Differentiate both side of the Poisson's equation $\nabla\cdot(\epsilon\nabla\phi)=\sum_{l=1}^M q^lc^l+\rho^f$ with regard to time $t$, then employ the identity $\boldsymbol{D}/\epsilon=-\nabla\phi$ and Equation \ref{NP_1} to obtain that $\frac{\partial \boldsymbol{D}}{\partial t}+\sum_{l=1}^M {q^l}\boldsymbol{J}^l$ is divergence-free. Representing this quantity with $\boldsymbol{\Theta}$, Equations \ref{A_1} and \ref{Theta_div_free} can be reached \cite{ANP000_1}.}

In \cite{ANP000_2} a standard while effective numerical scheme is proposed. {The main idea of \cite{ANP000_2} is to discretise Equations \ref{NP_1} and \ref{A_1} with the flux $\mathbf{J}^l$ reformulated following \cite{NP003}, then solve forward in time. During this process the deviation from Equation \ref{D_curl_free} is fixed by local curl-free relaxation algorithm.} At time $t=0$, the initial electric {potential $\phi^0$} is obtained with the {Poisson's equation:}

\begin{equation}
	{-\nabla\cdot(\epsilon\nabla\phi^0)=\sum_{l=1}^M q^lc^{l, \text{init}}+\rho^f,}\label{Poisson}
\end{equation}

where $c^{l, \text{init}}$ is the initial distribution of the $l$-th ions, and $\rho^f$ is the fixed charge. {Initial electric displacement $\boldsymbol{D}^0$ is calculated with a central difference scheme: }

\begin{equation}
	{\frac{D^0_{i+\frac12,j}}{\epsilon_{i+\frac12,j}}=\frac{\phi^0_{i+1,j}-\phi^0_{i,j}}{\Delta x},}
	{\frac{D^0_{i,j+\frac12}}{\epsilon_{i,j+\frac12}}=\frac{\phi^0_{i,j+1}-\phi^0_{i,j}}{\Delta y}}
\end{equation}

Equation \ref{NP_1} is discretised with a central difference scheme

\begin{equation}
\frac{c^{l,n+1}-c^{l,n}}{\Delta t}=-\frac{J^{l,n}_{i+\frac12,j}-J^{l,n}_{i-\frac12,j}}{\Delta x}-\frac{J^{l,n}_{i,j+\frac12}-J^{l,n}_{i,j-\frac12}}{\Delta y}\label{c_time_discretisation}
\end{equation}

where $J^{l,n}_{i\pm\frac12,j}$ and $J^{l,n}_{i,j\pm\frac12}$ are numerical approximations of $x$- and $y$-entries of flux function $\mathbf{J}^l$ at the $n$-th time step

\begin{equation}
{J^{l,n}_{i+\frac12,j}=-\frac{1}{\Delta x}\left(B(\Delta xq^l\frac{D^n_{i+\frac12,j}}{\epsilon_{i+\frac12,j}})-B(-\Delta xq^l\frac{D^n_{i+\frac12,j}}{\epsilon_{i+\frac12,j}})\right),\label{numerical_flux}}
\end{equation}

{where $B(x)=\frac{x}{e^x-1}$. Similarly for others.}

For $\boldsymbol{D}$, an intermediate value (denoted by $\boldsymbol{D}^*$) is first calculated {with} the discretised {Maxwell--}Amp\`{e}re equation{s} \ref{A_1} {and} \ref{NP_2_A_2}:

\begin{equation}
	{\frac{\boldsymbol{D}^*-\boldsymbol{D}^n}{\Delta t}=-\sum_{l=1}^M q^l\boldsymbol{J}^{l,n}+\boldsymbol{\Theta}^n}\label{Theta_and_D},
\end{equation}

{where $\boldsymbol{\Theta}^n=\frac{\boldsymbol{D}^n-\boldsymbol{D}^{n-1}}{\Delta t}+\sum_{l=1}^M \frac{q^l}{2\kappa^2}\boldsymbol{J}^{l,n-1}$ has been proven as the best approximation for $\boldsymbol{\Theta}^n$ through massive numerical tests \cite{ANP000_1}}. {Because} the value of $\boldsymbol{\Theta}$ is an approximation, the value of $\boldsymbol{D}^*$ is not the true value for the next time step. {This may violate Equation \ref{D_curl_free}}. Thus, starting from $\boldsymbol{D}^*$, an iterative constrained optimisation algorithm is performed to restore the correct value {of} $\boldsymbol{D}^{n+1}$:

\begin{equation}
{\boldsymbol{D}^{n+1}=\arg\min_{\boldsymbol{D}}\int_{\Omega}\frac{\Vert\boldsymbol{D}\Vert^2}{\epsilon}\text{ s.t. }\nabla\cdot\boldsymbol{D}=\sum_{l=1}^M q^lc^{l, n+1}+\rho^f.}\label{Local_curl_free_relaxation}
\end{equation}

{The minimiser of Equation \ref{Local_curl_free_relaxation} satisfies $\boldsymbol{D}+\epsilon\nabla\phi=0$ for some $\phi$ by Lagrange multiplier and calculus of variation \cite{NP021}, which means Equation \ref{D_curl_free} is satisfied.}

\section{Deep neural network and PINN}\label{sec:PINN}

Consider the PDE

\begin{equation}\label{PDE}
\begin{aligned}
	&\mathcal{N}u(\boldsymbol{x},t)=f(\boldsymbol{x},t),\boldsymbol{x}\in\Omega\subset\mathbb{R}^n,t\in[0,T]\\
	&\mathcal{B}u(\boldsymbol{x},t)=g(\boldsymbol{x},t),\boldsymbol{x}\in\partial\Omega,t\in[0,T]\\
	&u(\boldsymbol{x},0)=u_0(\boldsymbol{x}),\boldsymbol{x}\in\Omega
\end{aligned}
\end{equation}

where $u$ is the unknown solution function, $\mathcal{N}$ is a (probably nonlinear) differential operator, $f$ is a known inhomogeneous term, $g$ is a known function defining the {solution values} at the boundary $\partial\Omega$, $u_0$ is a known initial condition. The analytic solution is usually difficult to find even for {a} seemingly simple form of $\mathcal{N}${;} thus{,} the numerical methods for PDEs are then pursued. PINN {\cite{DLDE001}} has gained much attention in the community of scientific computing due to its portability and promising accuracy. The essence of PINN is the {utilisation} of a deep neural network as the approximator to fit the unknown solution function, similar to the polynomials {employed} in finite element methods. The basic building block of a deep neural network is called the layer, usually composed of an affine transformation and a nonlinear activation function. Mathematically, a layer $l_i:\mathbb{R}^{n_i}\rightarrow\mathbb{R}^{n_{i+1}}$, $i=1,2,...,L$ can be written as

\begin{equation}
l_i(\boldsymbol{x})=\sigma_i(\boldsymbol{W}_i\boldsymbol{x}+\boldsymbol{b}_i)\text{, }i=1,2,...,L,
\end{equation}

where $\boldsymbol{W}_i\in\mathbb{R}^{n_{i+1}\times n_i}$, $\boldsymbol{b}_i\in\mathbb{R}^{n_{i+1}}$, and $\sigma_i:\mathbb{R}^{n_{i+1}}\rightarrow\mathbb{R}^{n_{i+1}}$ is a nonlinear (usually element-wise) activation function. The whole deep neural network can be written as

\begin{equation}
\tilde{u}(\boldsymbol{x};\boldsymbol{\theta})=l_L\circ l_{L-1}\circ\dots\circ l_1(\boldsymbol{x}),
\end{equation}

where $\circ$ denotes the function composition, and $\boldsymbol{\theta}:=\{\boldsymbol{W}_i,\boldsymbol{b}_i\}_{i=1}^L$ is the collection of the network parameters. The variable $\boldsymbol{x}$ is the input to the neural network. In the settings of evolution equations, the input includes both the spatial variable and time, which {can be} packed{ }to form one input vector {(the continuous-time model)}. In the {remainder} of this section, we{ }abuse the variable{ }$\boldsymbol{x}$ to solely denote the spatial variable, thus the input to the network should be written as $(\boldsymbol{x},t)$.

To train the neural network means to optimise some loss function with regard to the network parameter $\{\boldsymbol{W}_i,\boldsymbol{b}_i\}_{i=1}^L$ {utilising some} optimisers{, such as} Stochastic Gradient Descent \cite{AWS04} or ADAM \cite{EEDL09}. In {the} PINN {\cite{DLDE001}}, the loss {function} contains the information of the governing equation {(see Equation \ref{PDE})} at collocation points $\{(\boldsymbol{x}_{i},t_{i})\}_{i=1}^{N_{eq}}\subset\Omega\times[0,T]$ {and paired data (usually obtained from initial and boundary conditions) $\{(\boldsymbol{x}_{D,i},t_{D,i}), u_{D,i}\}_{i=1}^N$:}

\begin{equation}\label{PINN_loss}
{\mathcal{L}=\frac{1}{N_{eq}}\sum_{i=1}^{N_{eq}} (\mathcal{N}\tilde{u}(\boldsymbol{x}_{i},t_{i})-f(\boldsymbol{x}_{i},t_{i}))^2+\frac{1}{N}\sum_{i=1}^{N} (\tilde{u}(\boldsymbol{x}_{D,i},t_{D,i})-u_{D,i})^2}
\end{equation}

The differentiation operation{s} appearing in $\mathcal{N}$ and $\mathcal{B}$ can be handled {with} AD \cite{ANP069}.

\section{Hybrid PINN}\label{sec:hybrid_PINN}

As remarkable as the original numerical method is, some improvements {can still be achieved}. {Firstly, the relaxation algorithm to correct the intermediate $\boldsymbol{D}^*$ requires changes in one cell to be compensated {for} in another cell{;} thus, it is inapplicable in {the} 1-dimensional case. Section \ref{sec:NE_1d} demonstrates that the same formula does not work for the same system in {the} 1-dimensional case. Therefore, a method to obtain an accurate intermediate value for the electric displacement is required. Secondly, it is far from obvious how to incorporate boundary condition involving $\boldsymbol{D}$ in the original conventional numerical method. For example, a common boundary condition for this system is the Neumann boundary condition} 

\begin{equation}
{\frac{\partial\phi}{\partial\boldsymbol{n}}(\boldsymbol{x})=g(\boldsymbol{x}),\boldsymbol{x}\in\partial\Omega,}
\end{equation}

{where $\boldsymbol{n}$ is the outward normal vector and $g$ is known. From the identity $\boldsymbol{D}/\epsilon=-\nabla\phi$, this boundary condition can be expressed in terms of $\boldsymbol{D}$:}

\begin{equation}
 {-\frac{\boldsymbol{D}}{\epsilon}\cdot\boldsymbol{n}=g.}\label{Neumann_BC}
 \end{equation}

 {Thirdly, the variable $\boldsymbol{\Theta}$, coming with the introduction of the novel Maxwell--Amp\`{e}re equation (Equations \ref{A_1} and \ref{NP_2_A_2}), brings an additional degree of freedom. There is no theoretical result for the value that should be chosen, except that it is divergence-free (see Equation \ref{Theta}). In \cite{ANP000_1, ANP000_2}, the value is selected via heuristics. They provide three options for $\boldsymbol{\Theta}$, and in their numerical tests the second one $\boldsymbol{\Theta}^n=\frac{\boldsymbol{D}^n-\boldsymbol{D}^{n-1}}{\Delta t}+\sum_{l=1}^M q^l\boldsymbol{J}^{l,n-1}$ is proven to be the best. It is curious because the third choice $\boldsymbol{\Theta}^n=\frac32(\frac{\boldsymbol{D}^n-\boldsymbol{D}^{n-1}}{\Delta t}+\sum_{l=1}^M q^l\boldsymbol{J}^{l,n-1})-\frac12(\frac{\boldsymbol{D}^{n-1}-\boldsymbol{D}^{n-2}}{\Delta t}+\sum_{l=1}^M q^l\boldsymbol{J}^{l,n-2})$ indeed gives higher order of accuracy. This suggests that figuring out principled guidance for approximations to this variable might be difficult. Clearly, considerable efforts have been invested in determining a delicate formula which will work practically. To make things worse, this formula is not universal, meaning that it merely applies to this particular system, and might be completely useless when the equations are further modified owing to advances in physics.}

To tackle the issues, we propose to use a neural network to calculate {the} values of the variable $\boldsymbol\Theta$. {This section presents} our proposed hybrid scheme for solving {M}ANP equations, analyse{s} its conservative properties, and describe{s} situations in 1-dimensional problems in details. 

\subsection{Approximation to $\boldsymbol{\Theta}$ through neural network}

{Unlike} vanilla PINN, for $\boldsymbol{\Theta}${,} we have neither paired data nor explicit governing physical laws{; hence,} a {novel approach} to design the loss function is required. {According to E}quation \ref{Theta_and_D}{, }$\boldsymbol{\Theta}$ and $\boldsymbol{D}$ are closely related{;} thus{,} we propose {controlling} the value of $\boldsymbol{D}${ }to indirectly guide the training {of} the neural network. For $\boldsymbol{D}$, there are governing {E}quations \ref{A_1}{and} \ref{NP_2_A_2}{. H}owever, {because} many variables are involved,{ }numerical error{s} from {these} variables may {accumulate}. {In addition,} it {slows} down the training {because} the loss value and its partial derivatives with {regard} to every network {parameter} must be calculated in each training loop. To improve efficiency, we {apply} Equation \ref{D_curl_free} to design the loss function. The straightforward formulation is {as follows}:

\begin{equation}
	\mathcal{L}_{PI}{(t)}=\frac{1}{N}\sum_{i=1}^N \left(\frac{\partial}{\partial x}\frac{D^{(2)}}{\epsilon}(x_i,y_i, t)-\frac{\partial}{\partial y}\frac{D^{(1)}}{\epsilon}(x_i,y_i, t)\right)^2,\label{L_PI_original}
\end{equation}

or an alternative loss design {adopting} Lagrange multiplier to eliminate the derivative term \cite{NP021}{ }:

\begin{equation}
	{\mathcal{L}^{\prime}_{PI}{(t)}}=\frac{1}{N}\sum_{i=1}^N \left((\frac{D^{(1)}}{\epsilon}(x_i,y_i, t))^2+(\frac{D^{(2)}}{\epsilon}(x_i,y_i, t))^2\right)\label{L_PI}
\end{equation}

where $D^{(1)}$ and $D^{(2)}$ are $x$- and $y$-entry of the vector field $\boldsymbol{D}$, respectively, $\{(x_i,y_i)\}_{i=1}^N$ are the spatial locations used in training, and $t$ is the current time step. Equation \ref{L_PI} does not contain derivative terms{; therefore, it is} more {numerically stable}. In the numerical experiments, we find that Equation \ref{L_PI} {gives} better results.
\\ \hspace*{\fill}

\begin{remark}
To be rigorous, the minimiser of Equation \ref{L_PI} is not the solution satisfying $\frac{\boldsymbol{D}}{\epsilon}=-\nabla\phi$ and the ${\nabla\cdot\boldsymbol{D}=\sum_{l=1}^M q^lc^{l}+\rho^f}$. The true solution is the constrained minimiser under the condition {${\nabla\cdot\boldsymbol{D}=\sum_{l=1}^M q^lc^{l}+\rho^f}$}. From Equation \ref{Theta_and_D}, we can see that given {$\boldsymbol{\Theta}$ is divergence-free {in the discrete sense}, and the current $\boldsymbol{D}$ satisfies the this constraint}, the updated $\boldsymbol{D}^*$ will {do so}. It will be shown that in our proposed method, the divergence-free condition can be automatically and analytically guaranteed{;} thus{,} we can safely discard the constraint, and{ }perform unconstrained optimisation, for which many techniques exist in deep learning contexts.
\end{remark}

{As mentioned before, $\boldsymbol{\Theta}$ and $\boldsymbol{D}$ are closely related, we apply the numerical value of the previous time step as the input of the neural network \cite{ANP071, ANP074}.} Inspired by Theorem 3.1 in \cite{hANP001}, {we} train it to fit a hidden scalar field $\tilde{u}(x,y,t)$ then obtain the entries $(\Theta^{(1)}, \Theta^{(2)})$ of the target $\boldsymbol{\Theta}$ through

\begin{equation}
\begin{aligned}
\Theta^{(1)}(x_i,y_i,t)&={\frac{\tilde{u}(x_i,y_{i+1},t)-\tilde{u}(x_i,y_i,t)}{y_{i+1}-y_i}}\\
\Theta^{(2)}(x_i,y_i,t)&={-\frac{\tilde{u}(x_{i+1},y_{i},t)-\tilde{u}(x_i,y_i,t)}{x_{i+1}-x_i}.}\\
\end{aligned}
\end{equation}

{Using} this formulation, the divergence-free condition is {satisfied} automatically {in the discrete sense.} Thus, we only need to minimise {$\mathcal{L}^{\prime}_{PI}$} (Equation \ref{L_PI}) during training the neural network. 

{Because deep learning is an optimisation-based numerical approximation method, it is highly flexible to implement almost any boundary condition, as a penalty term. Take the Neumann boundary condition (Equation \ref{Neumann_BC}) as an example, the penalty term can be written as:}

\begin{equation}
{\mathcal{L}_{BC}(t)=\frac{1}{N_{BC}}\sum_{i=1}^{N_{BC}}\left(-\frac{\boldsymbol{D}(x_i, y_i, t)}{\epsilon}\cdot\boldsymbol{n}+g(x_i, y_i, t)\right)^2,}\label{L_BC}
\end{equation}

{where $\{x_i, y_i\}_{i=1}^{N_{BC}}$ are the grid points at the boundary. In practice, because the values of $\Theta^{(1)}$ and $\Theta^{(2)}$ should not variate too much in a single discrete cell, to avoid overfitting and instability, we add a regularisation term: }

\begin{equation}
\begin{aligned}
{\mathcal{R}(t)=}&{\int_\Omega\Vert\nabla\Theta^{(1)}(x,y,t)\Vert_2^2+\Vert\nabla\Theta^{(2)}(x,y,t)\Vert_2^2dxdy.}\\
{\approx}&{\Delta_x\Delta_y\sum_{i=1}^N\left(\Vert\nabla\Theta^{(1)}(x_i,y_i,t)\Vert_2^2+\Vert\nabla\Theta^{(2)}(x_i,y_i,t)\Vert_2^2\right).}\label{Regularisation}
\end{aligned}
\end{equation}

{The final loss function for the neural network training is: }

\begin{equation}
{\mathcal{L}(t)=\mathcal{L}^{\prime}_{PI}(t)+\lambda_{BC}\mathcal{L}_{BC}(t)+\lambda_R\mathcal{R}(t),}
\end{equation}

{where $\lambda_{BC}$ and $\lambda_{R}$ are hyper-parameters.}

{During training, at each time step $t_i$, we update the network parameters until the loss function value $\mathcal{L}(t_i)$ is small enough or a preset training loop threshold is reached.}

\subsection{Analysis of conservative properties}

{This section analyses the conservation properties of the proposed method: conservation of total mass and preservation of ion concentration positivity. Although a neural network replaces the concrete formula for $\boldsymbol{\Theta}$, which introduces stochasticity (e.g., during initialisation and training), the conservation properties are retained, which is different from vanilla pure end-to-end PINN. Such conservation properties are generally difficult to ensure in vanilla PINN. In situations where the network is perfectly trained and converges to the true solution exactly we are certain that these properties are possessed by the numerical solutions given by the neural network. However, in most cases, deep learning practitioners do not have this guarantee. In the proposed method, these conservation properties still hold even if approximation errors are present (e.g. owing to insufficient training).}

{ }Theorem 3.1 of \cite{ANP000_2} {establishes} the mass conversation property through summing up both sides of {E}quation \ref{c_time_discretisation} over the {entire} spatial domain. In our method, the values of $\boldsymbol{J}$ might be different, {because} the values given by the neural network at different training stages affect the value of $\boldsymbol{D}$ and{,} consequently the values of $\boldsymbol{J}$. However, the mass {is} still {conserved because} all these {$\boldsymbol{J}$ cancel} provided {a periodic} boundary condition.

The mass positivity of the original algorithm is proven in Theorem 3.2 of \cite{ANP000_2}. They first flatten the {2D} array containing {the ion concentration values} at grid points to be{ Vector} $c^{n,l}$, then construct {Matrix} $\mathcal{L}$ such that $\mathcal{L}c^{n+1,l}=c^{n,l}$. The particular matrix {entry} values obtained in our hybrid method could be different, but {this matrix} remains strictly{ }diagonally dominant (thus invertible){,} with positive diagonal entries and negative off-diagonal entries. {According to} Gershgorin's circle theorem, all{ }eigenvalues have positive real parts{;} thus{, it} is an M-matrix. {According to} statement $F_{15}$ {in} \cite{hANP002}, its inverse has no negative entries{;} thus{,} this transformation {preserves} positivity.

\subsection{Generalisation to \texorpdfstring{{1-dimensional }}{1-dimensional }problems}\label{sec:gen}

{Generalising} this particular numerical scheme to {1-dimensional problems is not trivial}. The first problem is with the empirical equation for the variable $\boldsymbol{\Theta}$. Because this expression is obtained from extensive experiments in {a} 2-dimensional setting, {its transfer to 1-dimensional cases remains unclear}. {We} will {demonstrate} that in 1-dimensional problems, neither directly adapting the formula of $\boldsymbol{\Theta}$ in {the} 2-dimensional case, nor simply setting it {as} 0 will{ }work in numerical experiments. {The other reason accounting for the failure of the original algorithm in 1-dimensional problems} is that the formula is an approximation for correct $\boldsymbol{\Theta}$ values. {However, this} approximation {could be inaccurate}. {To fix possible error caused by this imperfect approximation of $\boldsymbol{\Theta}$,} in \cite{ANP000_1, ANP000_2}, the authors {apply} a {local curl-free} relaxation algorithm to correct it. This procedure is simple yet effective in 2-dimensional problems, {but in} 1-dimensional case{s}, this algorithm is {inapplicable}, because once we update the electric displacement in one cell, there is no cell to compensate {for} it to restore the {constraint in Equation \ref{Local_curl_free_relaxation}}. The proposed hybrid scheme provides an efficient means to adaptively calculate the proper values for $\boldsymbol{\Theta}$ in different settings {via} neural network training,{ }an automatic and user-transparent optimisation process.

{Here we describe the proposed method in a 1-dimensional case.} The PNP system in {a} 1-dimensional space is:

\begin{numcases}{}
	\frac{\partial c^l}{\partial t}= -\frac{\partial J^l}{\partial x}\text{, }l=1,2,...,M\label{NP_1_1d}\\ 
	J^l=-(\frac{\partial c^l}{\partial x}+q^lc^l\frac{\partial \phi}{\partial x})\text{, }l=1,2,...,M\label{NP_2_A_2_1d}\\
	\epsilon_0^2\frac{\partial^2 \phi}{\partial x^2}=-\rho^f-\sum_{l=1}^M q^lc^l\label{Poisson_1d},
\end{numcases}

where $\epsilon_0$ is a constant. Following the derivation of \cite{ANP000_1}, {we take the time-derivative of both sides of Equation \ref{Poisson_1d}, plugging in Equation \ref{NP_1_1d}, and using $\Theta$ to represent $\epsilon_0^2\frac{\partial}{\partial t}\frac{\partial \phi}{\partial x}-\sum_{l=1}^M q^lJ^l$}, the {Maxwell--}Amp\`{e}re equation is:

\begin{numcases}{}
	\epsilon_0^2\frac{\partial}{\partial t}\frac{\partial \phi}{\partial x}=\Theta+\sum_{l=1}^M q^lJ^l\\
	\frac{\partial \Theta}{\partial x}=0.
\end{numcases}

Following \cite{NP020}, we consider the spatial domain $[-1, 1]$. For ionic flux $J^l$'s, no flux boundary condition is applied: 

\begin{equation}
J^l(-1)=J^l(1)=0\text{, }l=1,2,...,M,
\end{equation}

while a Robin-type boundary condition is applied for electric potential $\phi$:

\begin{numcases}{}
\phi(1)+\eta\frac{\partial \phi}{\partial x}(1)=\phi_0(1)\label{right_RBC}\\
\phi(-1)-\eta\frac{\partial \phi}{\partial x}(-1)=\phi_0(-1)\label{left_RBC},
\end{numcases}

where $\eta$, $\phi_0(1)$ and $\phi_0(-1)$ are constants.

The original numerical scheme solves the {M}ANP equations for $c^l$, $J^l$ and $\frac{\partial \phi}{\partial x}$. Similar to 2-dimensional case, if the value of $\Theta$ is not correctly selected, the resulting $\frac{\partial \phi}{\partial x}$ may not corresponding to an existing $\phi$, because the Robin-type boundary conditions for $\phi$ at {the} right and left end{s} could contradict. Thus, a loss function {reflecting} such contradiction should be designed. Consider{ing} the difference {between} both sides of Equation{s} \ref{right_RBC} {and \ref{left_RBC}}, we obtain:

\begin{equation}
\phi(1)-\phi(-1)+\eta(\frac{\partial \phi}{\partial x}(1)+\frac{\partial \phi}{\partial x}(-1))-(\phi_0(1)-\phi_0(-1))=0.
\end{equation}

After $\frac{\partial \phi}{\partial x}$ is calculated, we use the values of $\frac{\partial \phi}{\partial x}$ and rectangular rules of numerical integration to approximate $\phi(1)-\phi(-1)$:

\begin{equation}
\phi(1)-\phi(-1)=\int_{-1}^1\frac{\partial \phi}{\partial x}(x)dx\approx\sum_{i=1}^m\frac{\partial \phi}{\partial x}(x_i),
\end{equation}

where $-1=x_1<x_2<\cdots<x_m$ are {the} grid points used for numerical integration, and $m$ is the total number of{ }grid points.

We can write out the loss function as a squared error:

\begin{equation}
\mathcal{L}_{1d}=\left(\sum_{i=1}^m\frac{\partial \phi}{\partial x}(x_i)+\eta(\frac{\partial \phi}{\partial x}(1)+\frac{\partial \phi}{\partial x}(-1))-(\phi_0(1)-\phi_0(-1))\right)^2.
\end{equation}

{Considering that this loss function itself takes the boundary condition into account, and the output of the network is a scalar, we simply discard the loss for the boundary (Equation \ref{L_BC}) and the regularisation term (Equation \ref{Regularisation}).}

\section{Numerical experiments}\label{sec:NE}

\subsection{Steady\texorpdfstring{{-}}{-}state solution in 1-dimensional space}\label{sec:NE_1d}

{This section presents} the numerical results of our proposed method{ }and validate it through its convergence to the steady{-}state solution as $t\rightarrow\infty$ {following \cite{NP002}}. The steady{-}state solution is calculated {with} the method provided by {Lee et al} \cite{NP020}. We also compare with the experimental formulas for $\Theta$ from the original method (together with a variant){,} which have been proven effective in 2-dimensional case{s}:

\begin{equation}
\Theta^n=0\label{theta_0},
\end{equation}

\begin{equation}
\Theta^n=	\epsilon_0^2\frac{\frac{\partial \phi}{\partial x}^n-\frac{\partial \phi}{\partial x}^{n-1}}{\Delta t}-\sum_{l=1}^M q^lJ^{l, n-1}\label{2d_gen},
\end{equation}

\begin{equation}
\Theta^n=	\epsilon_0^2\frac{\frac{\partial \phi}{\partial x}^n-\frac{\partial \phi}{\partial x}^{n-1}}{\Delta t}-\sum_{l=1}^M q^lJ^{l, n}\label{2d_gen_2},
\end{equation}

where $\Delta t$ is the time step interval, and the superscript $n$ refers to the index for the time step.

The frame of the numerical scheme is a straightforward generalisation from the original 2-dimensional numerical method. Assume $\rho^f=0$ and $\epsilon_0=2^{-2}$, we first initialise $\phi^0$ and $\frac{\partial \phi}{\partial x}^0$ such that they satisfy the Poisson{'s} equation (Equation \ref{Poisson_1d}). {Numerical} values of $c^l$'s and $\frac{\partial \phi}{\partial x}$ are updated for each time step. After $\frac{\partial \phi}{\partial x}$ is solved, we solve for $\phi$ according to $\frac{\partial \phi}{\partial x}$ and the boundary condition at the left end {utilising} Euler forward ordinary differential equation solver.

\begin{figure*}
  \centering
  \subcaptionbox{}{\includegraphics[scale=0.2]{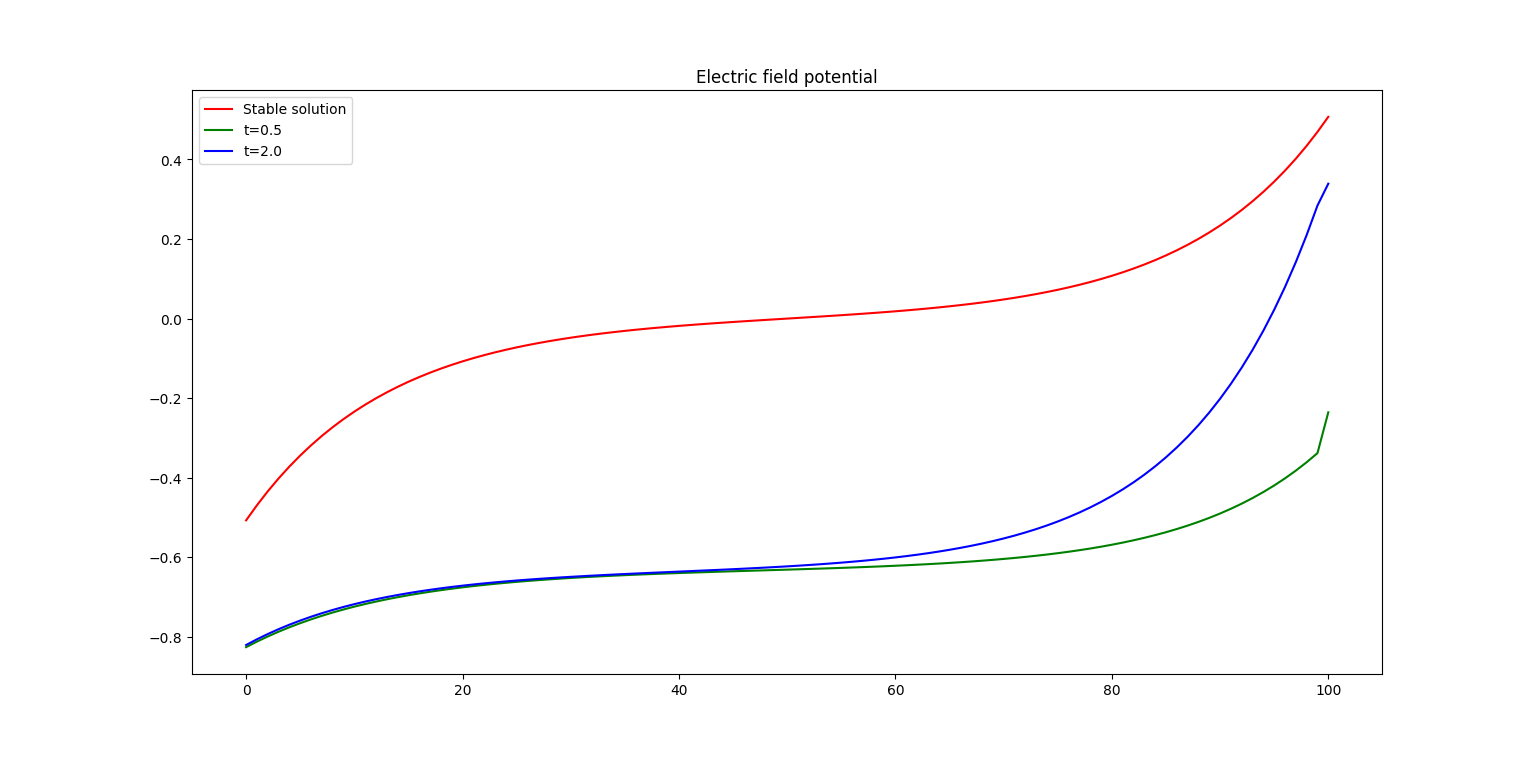}}
  \hfill
  \subcaptionbox{}{\includegraphics[scale=0.2]{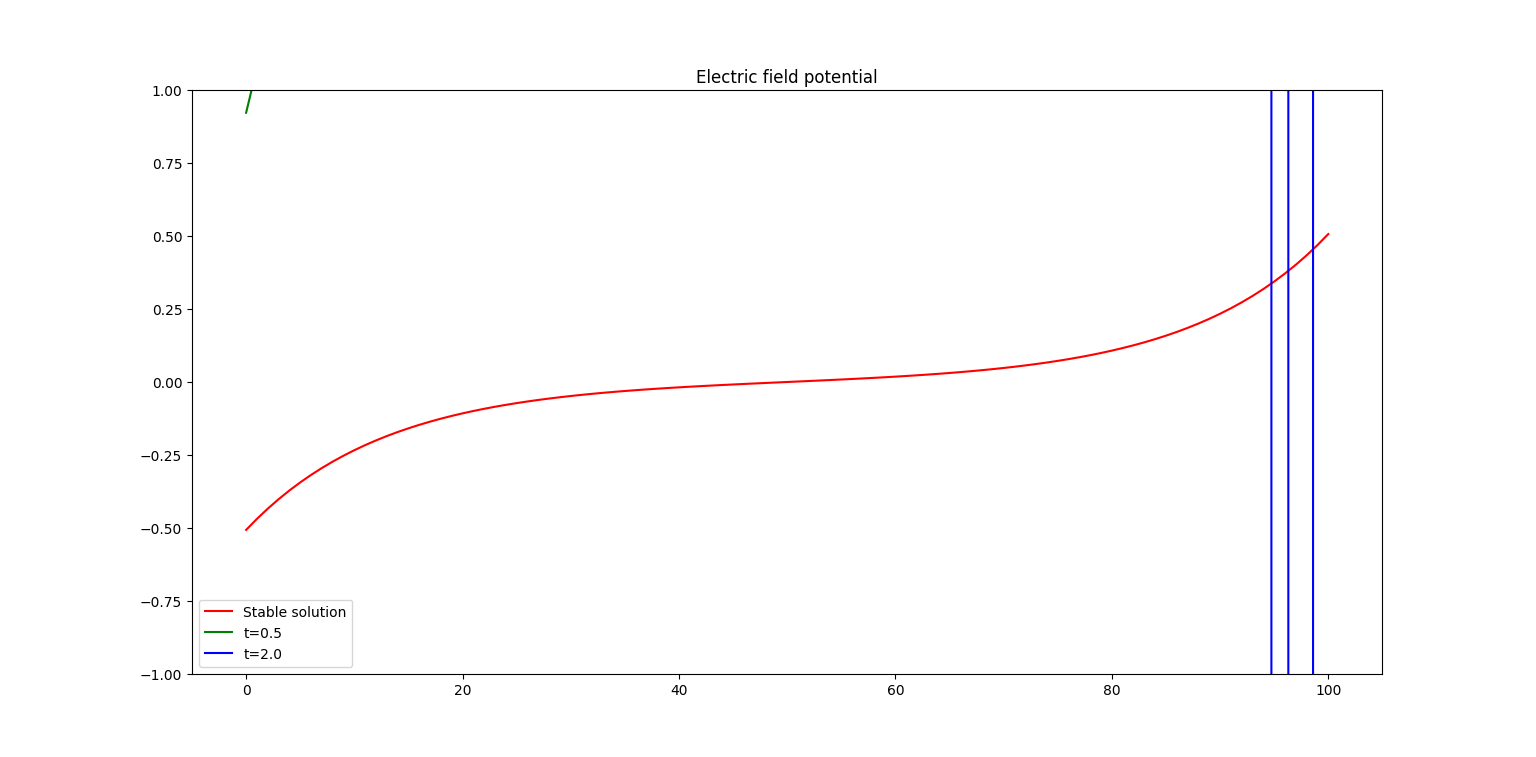}}
  \\
  \subcaptionbox{}{\includegraphics[scale=0.2]{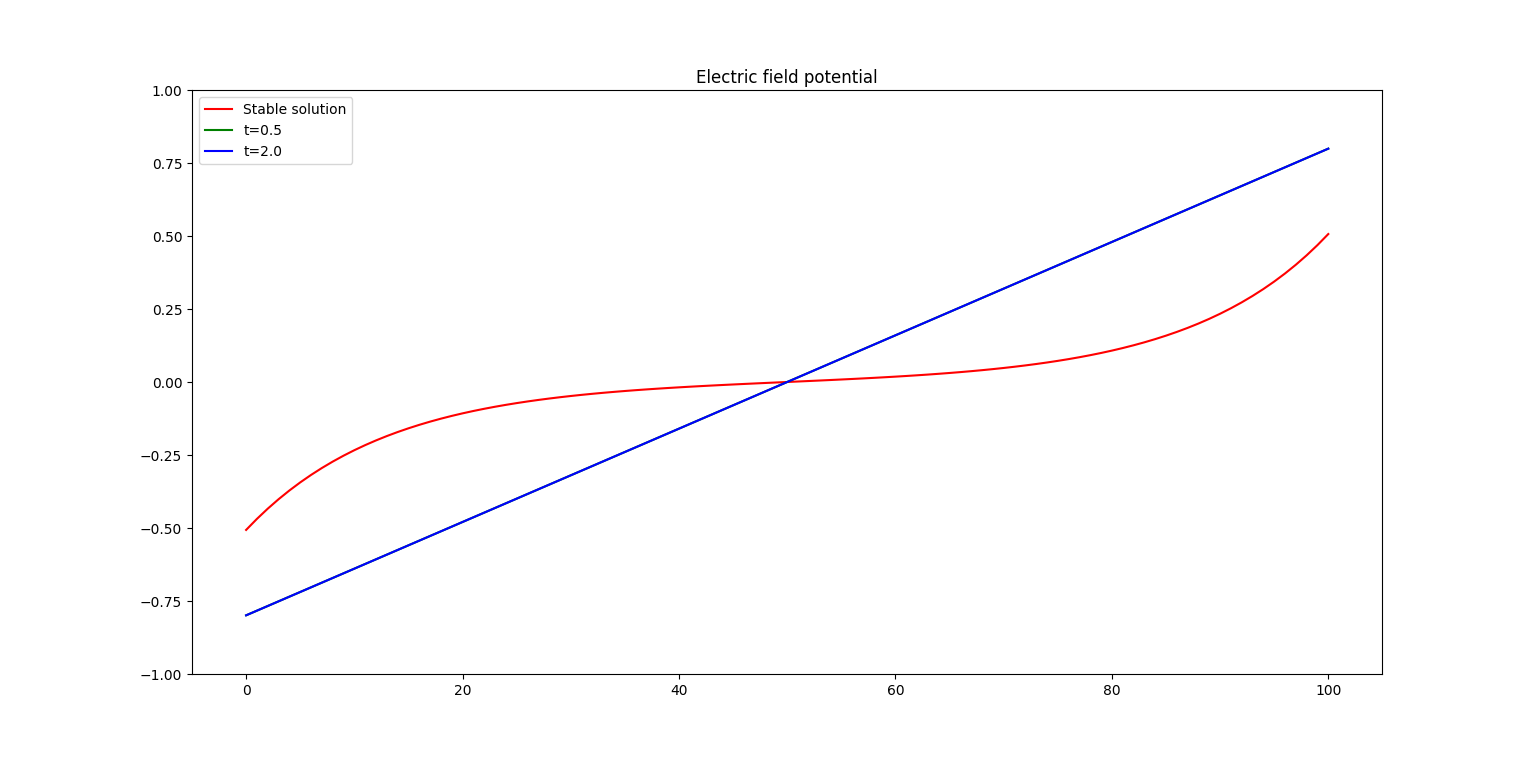}}
  \hfill
  \subcaptionbox{}{\includegraphics[scale=0.2]{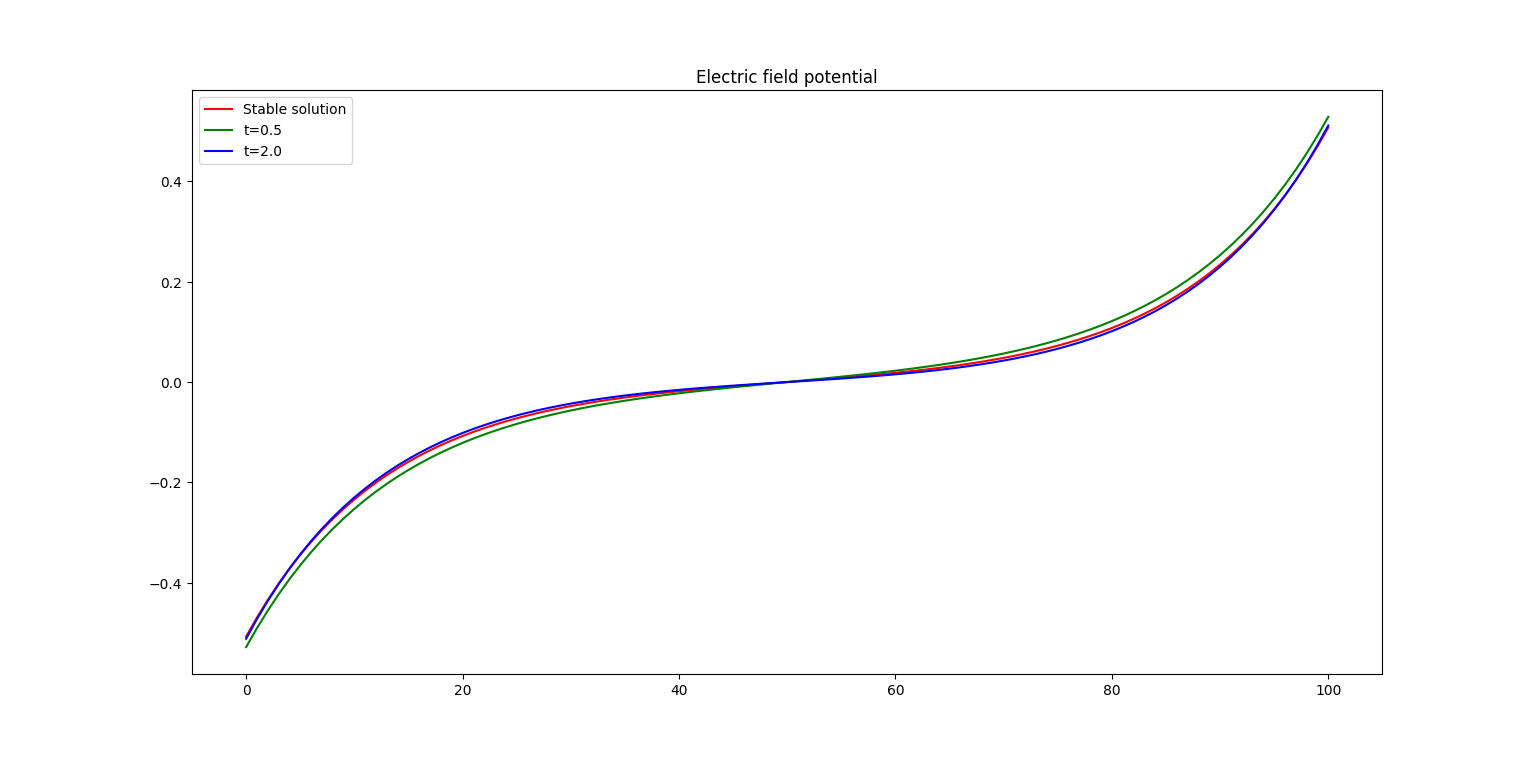}}
\caption{{E}lectric potential $\phi$ at different time {steps}. (a) $\Theta$ is set to be 0; (b) $\Theta$ is calculated {with} Equation \ref{2d_gen}; (c) $\Theta$ is calculated {with} Equation \ref{2d_gen_2}, the solution at $t=0.5$ and $t=2.0$ are overlapped{ }; (d) $\Theta$ is approximated by a neural network (ours)}
\label{fig:1d_all}
\end{figure*}

Figure \ref{fig:1d_all} {illustrates} the convergence behaviour of the numerical solutions {over time}. Simply setting $\Theta$ to be $0$ will not give the correct solution. The direct generalisation of the formula for $\Theta$ in 2-dimensional case{s} (Equation \ref{2d_gen}) exhibits numerical instability in {the} 1-dimensional problem. {When} a variant{ }(Equation \ref{2d_gen_2}) is applied, the numerical scheme fails to update the solution correctly {over time}. {The} solution{s} at different time {steps} are {identical}. {By contrast}, the solution calculated {with} the proposed method gradually converges to the steady{-}state solution. 

\begin{figure}
  \centering
  \includegraphics[scale=0.3]{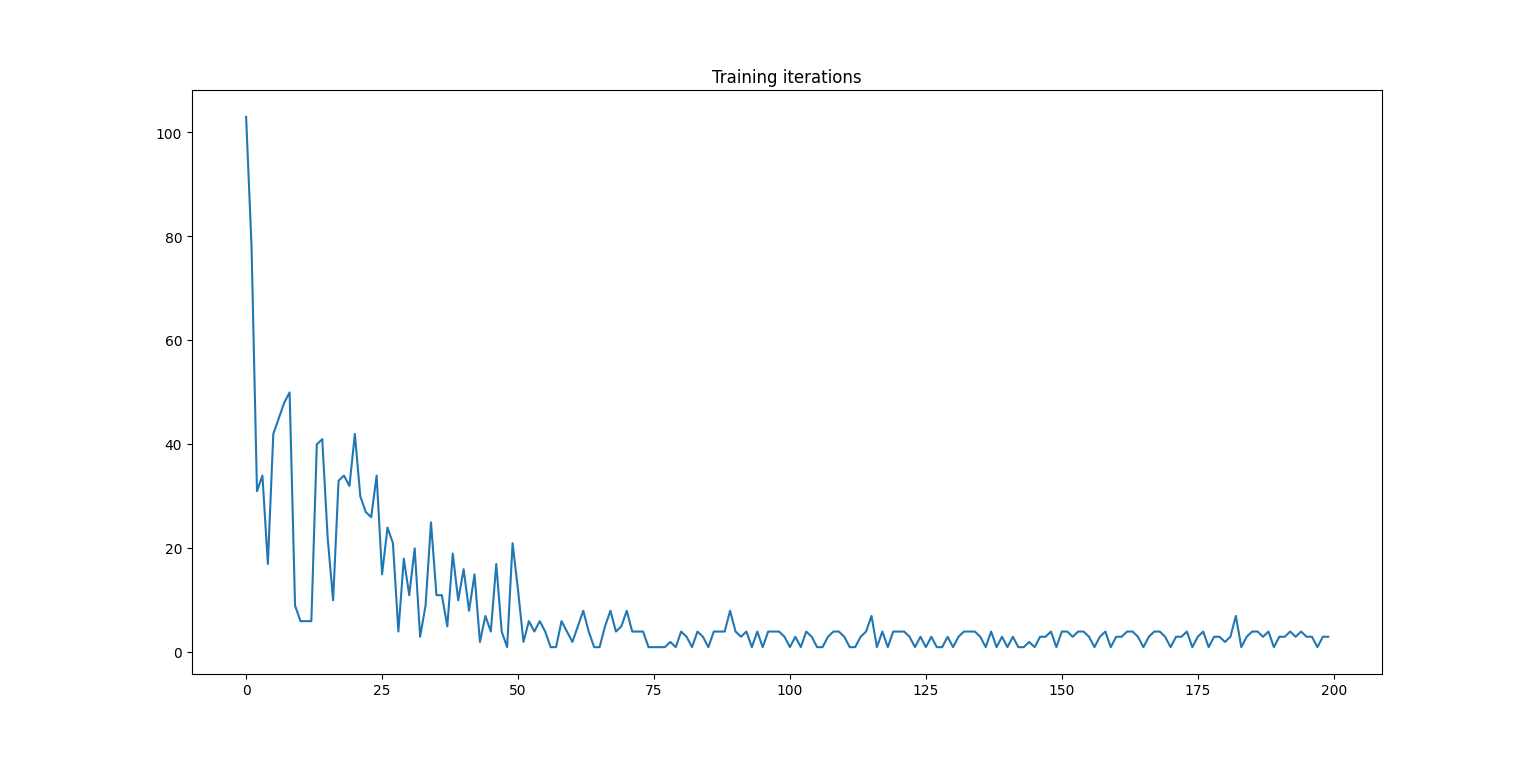}
  \caption{Training iterations required to reach a loss less than $10^{-8}$}
  \label{fig:training_iters}
\end{figure}

Figure \ref{fig:training_iters} shows the training iterations {that} the neural network {required} to optimise the loss value to{ }less than $10^{-8}$. {Although} the number of iterations is large at the beginning of {the} training, it drops quickly as the time step advances. 

There might be an argument that{,} in this simple problem, the correct optimiser of the loss function can be found analytically{ }without the assistance of a neural network. However, {this} depends {on} the discretisation methods {employed} for the PDE {and} the numerical integration schemes applied. {This may not be the case if a more sophisticated numerical method is adopted}. Furthermore, this approach {can be} immediately {applied to} another PDE system{ }in which a dummy variable providing an extra degree of freedom without a clear physical meaning{ }is present. {Researchers} need {not} conduct extensive numerical tests to {determine} a well-behaved {approximation formula. A} correction algorithm is {not required} to augment the numerical scheme{ }to avoid{ }error {accumulation}.

\subsection{\texorpdfstring{{Efficiency analysis}}{Efficiency analysis}}\label{sec:Efficiency}

{In the original work the local curl-free algorithm is performed cell-by-cell \cite{ANP000_1, ANP000_2}. This implementation is not suitable for parallelisation. We propose calculating the update step sizes for the all cells at once, organising them as matrices, and applying the updates to the whole domain. Although this trade a little accuracy for the efficiency, our proposed method still achieves higher accuracy (see Section \ref{sec:Analytic}), and meanwhile, the computation time is greatly reduced. Table \ref{tab:com_time} shows the running time of our method and the original method in \cite{ANP000_1, ANP000_2}, both computing $100$ time steps from $t=0$ to $t=0.5$. The data is measured on an Intel(R) Core(TM) i9-13900HX CPU (no CUDA \cite{AWS05} is utilised for the neural network). It can be observed that even our method includes the overhead of neural network training and inference, with the accelerated local curl-free algorithm, it is still faster, and the edge is larger when the grid is finer. Note that the proposed method is an integrated process including both formula discovery and numerically solving the system, but when running the original method, a well-tested formula to approximate $\Theta$ found in \cite{ANP000_1, ANP000_2} is already given.}

\begin{table}
 \caption{Running time}
  \centering
  \begin{tabular}{lll}
    \toprule
       Spatial resolution  & Ours     & Original \\
    \midrule
    $50\times50$ & 1m41.5s  & 1m49.7s     \\
    $100\times100$     & 6m55.2s & 10m46.7s      \\
    \bottomrule
  \end{tabular}
  \label{tab:com_time}
\end{table}

\subsection{\texorpdfstring{{Analytic tests}}{Analytic tests}}\label{sec:Analytic}
{Consider the domain $\Omega=[-1, 1]\times[-1, 1]$, and $t\in[0, 0.5]$. Let $q^1=+1$, $q^2=-1$, and $\epsilon=1$, we assume the exact solution of the system is:}

\begin{numcases}{}
	\phi_{true}(\boldsymbol{x},t)=\frac12 \Vert\boldsymbol{x}\Vert_2^2e^{-t}\\
	c^1_{true}=e^{-\phi}\\
	c^2_{true}=e^{\phi}\\
	\boldsymbol{D}_{true}=-\epsilon\nabla\phi_{true}.
\end{numcases}

{The system is then adapted to fit the exact solution above:}

\begin{numcases}{}
	{\frac{\partial c^l}{\partial t}= -\nabla\cdot\boldsymbol{J}^l+f_l\text{, }l=1,2}\\ 
	{\frac{\partial \boldsymbol{D}}{\partial t}=-\sum_{l=1}^2 q^l\boldsymbol{J}^l+\boldsymbol{\Theta}+\boldsymbol{h}}\\
	{\boldsymbol{J}^l=-(\nabla c^l-\frac{q^lc^l}{\epsilon}\boldsymbol{D})\text{, }l=1,2,}
\end{numcases}

{where $f_l$'s are calculated through the exact solution, and $\boldsymbol{h}$ is selected to be:}

\begin{equation}
{\boldsymbol{h}=}\left(
\begin{array}{l}
{x_1-x_2}\\
{x_2-x_1}
\end{array}
\right){e^{-t},}
\end{equation}

{so that $\boldsymbol{\Theta}$ is divergence-free (Equation \ref{Theta_div_free}). For the $c^l$'s, we apply periodic boundary condition, and for $\boldsymbol{D}$ we apply boundary condition in the Equation \ref{Neumann_BC}, which is the Neumann boundary condition for the electric potential $\phi$. The right-hand side of the Equation \ref{Neumann_BC} is obtained from the exact solutions. Suppose the solution given by the numerical algorithm are $c^l_{numerical}$'s and $\boldsymbol{D}_{numerical}$, the error of time $t$ is then calculated via:}

\begin{equation}
{E_{c^l}(t)=\sqrt{\frac{1}{N}\sum_{i=1}^N(c^l_{numerical}(x_i, y_i, t)-c^l_{true}(x_i, y_i, t))^2}}\label{c_error}
\end{equation}

\begin{equation}
{E_{\boldsymbol{D}}(t)=\frac{1}{N}\sum_{i=1}^N\Vert\boldsymbol{D}_{numerical}(x_i, y_i, t)-\boldsymbol{D}_{true}(x_i, y_i, t)\Vert_2.}\label{D_error}
\end{equation}

\begin{figure*}
\centering
	\subcaptionbox{}{\includegraphics[scale=0.1]{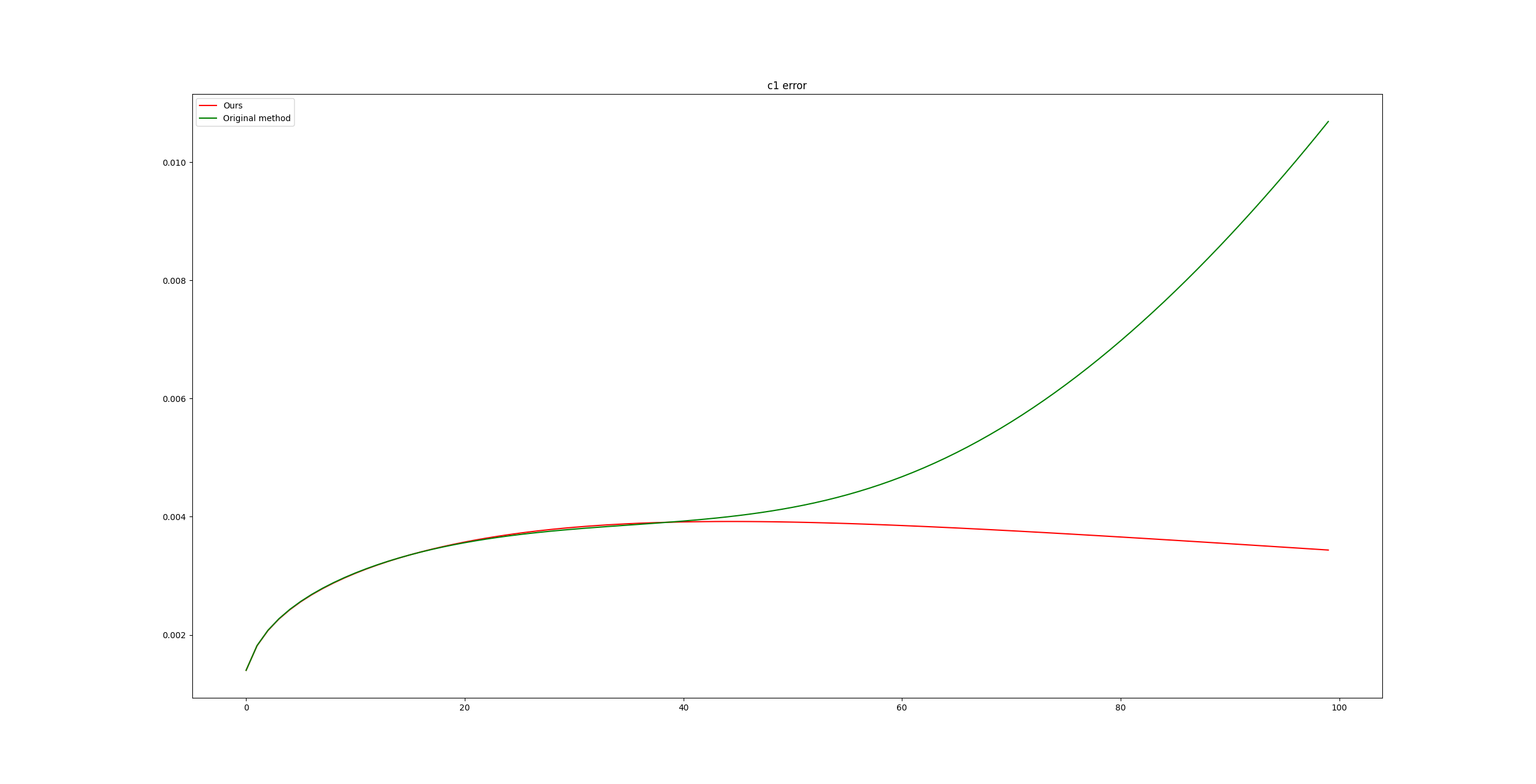}}
	\subcaptionbox{}{\includegraphics[scale=0.1]{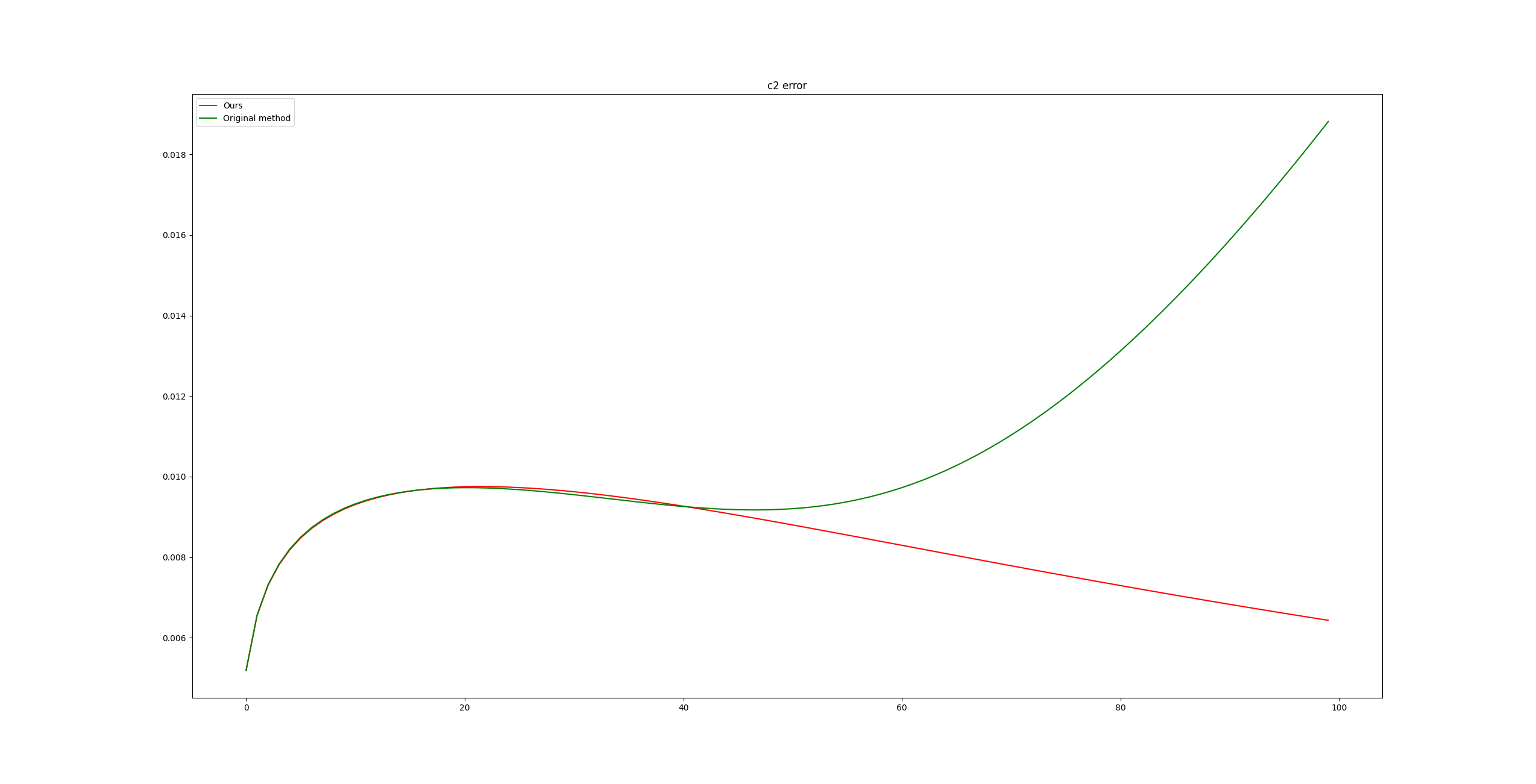}}
	\subcaptionbox{}{\includegraphics[scale=0.1]{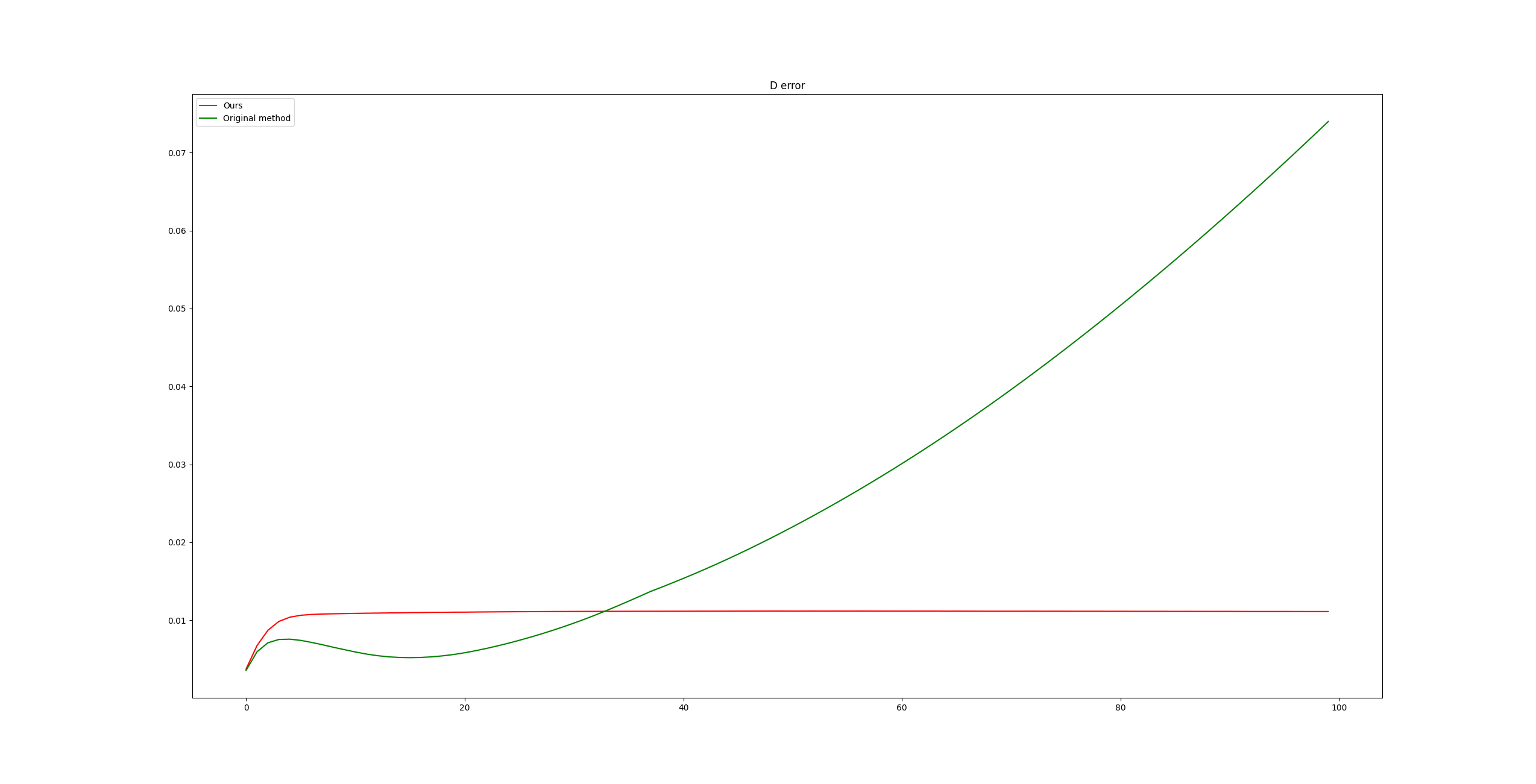}}
\centering
	\subcaptionbox{}{\includegraphics[scale=0.1]{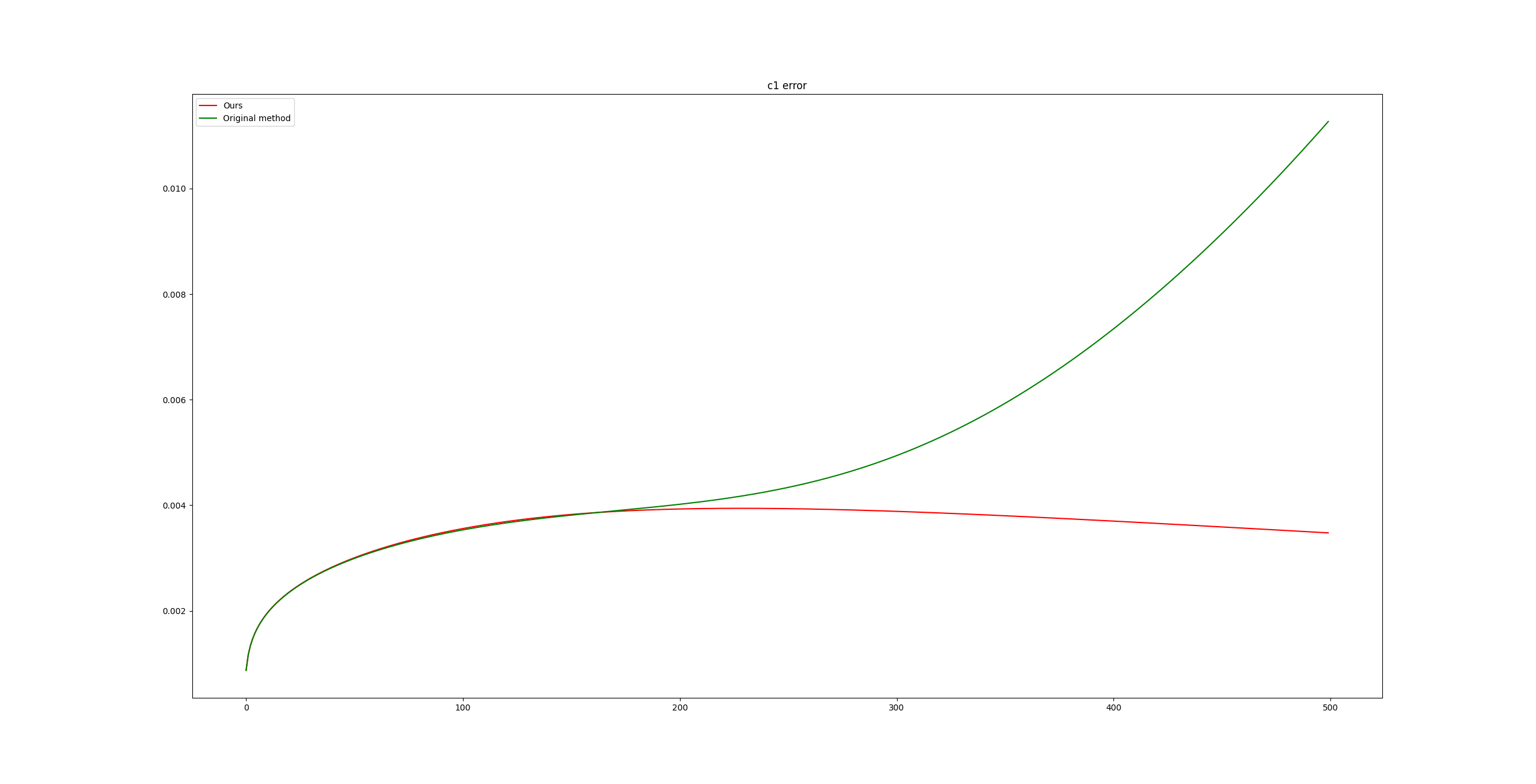}}
	\subcaptionbox{}{\includegraphics[scale=0.1]{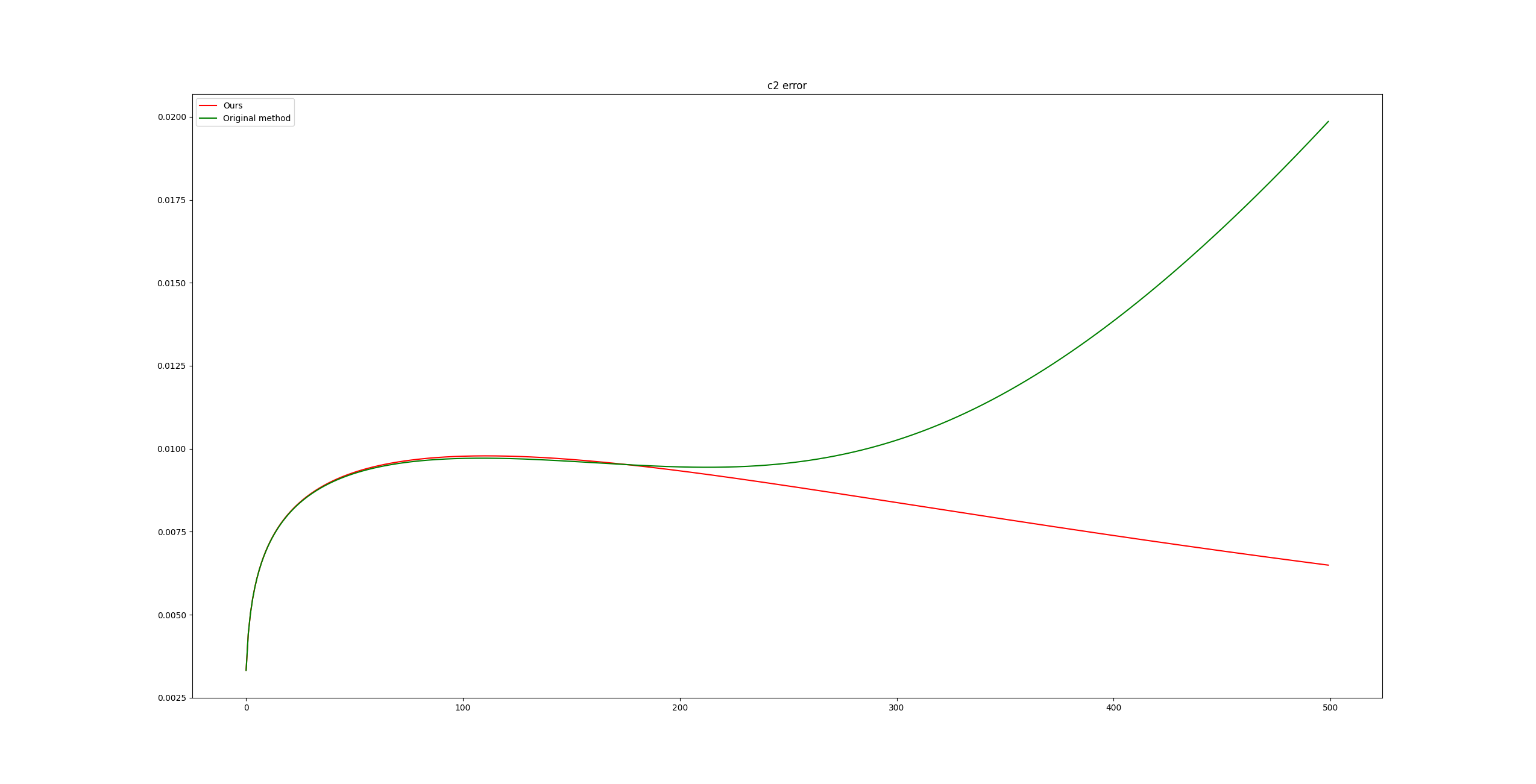}}
	\subcaptionbox{}{\includegraphics[scale=0.1]{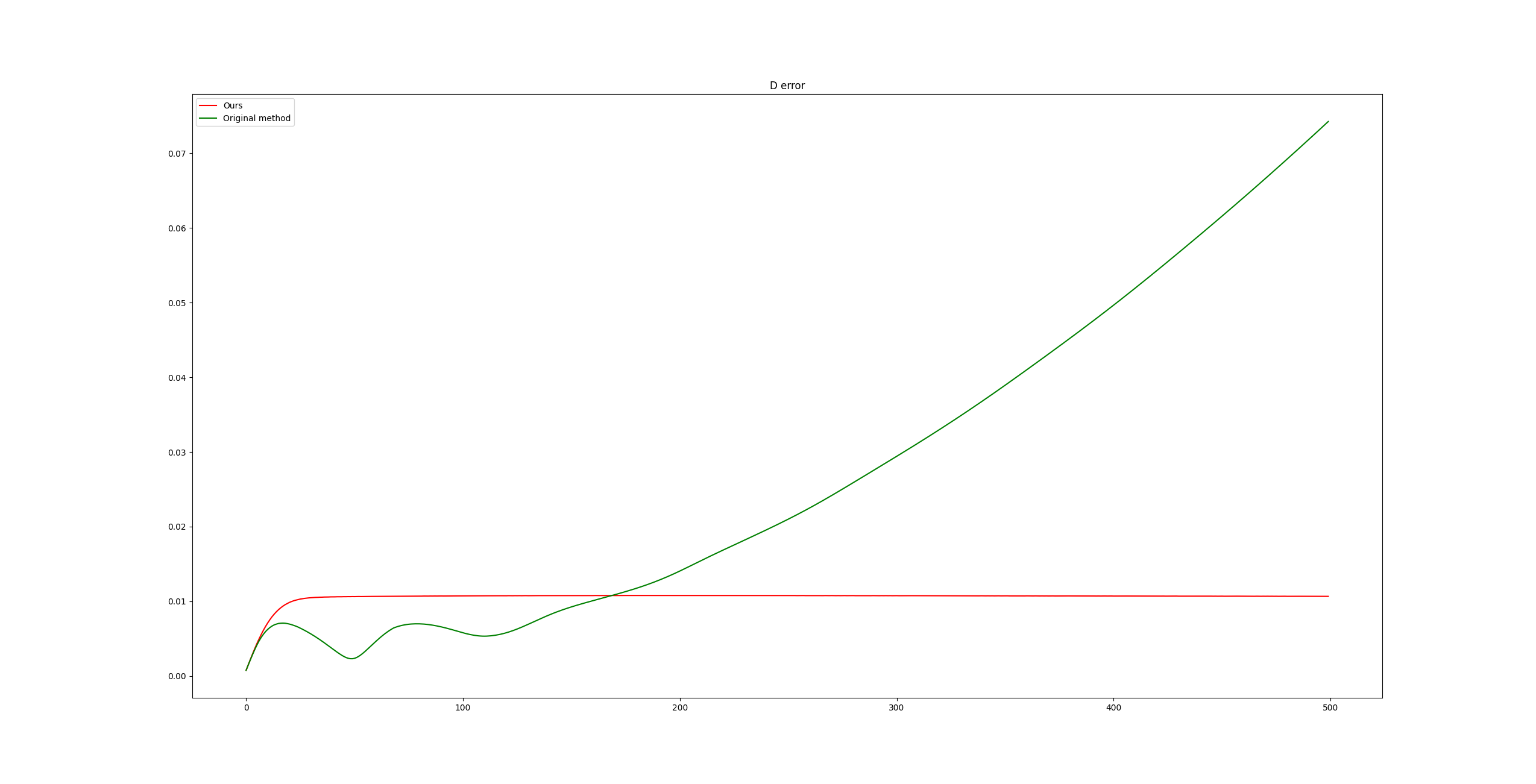}}
\centering
	\subcaptionbox{}{\includegraphics[scale=0.1]{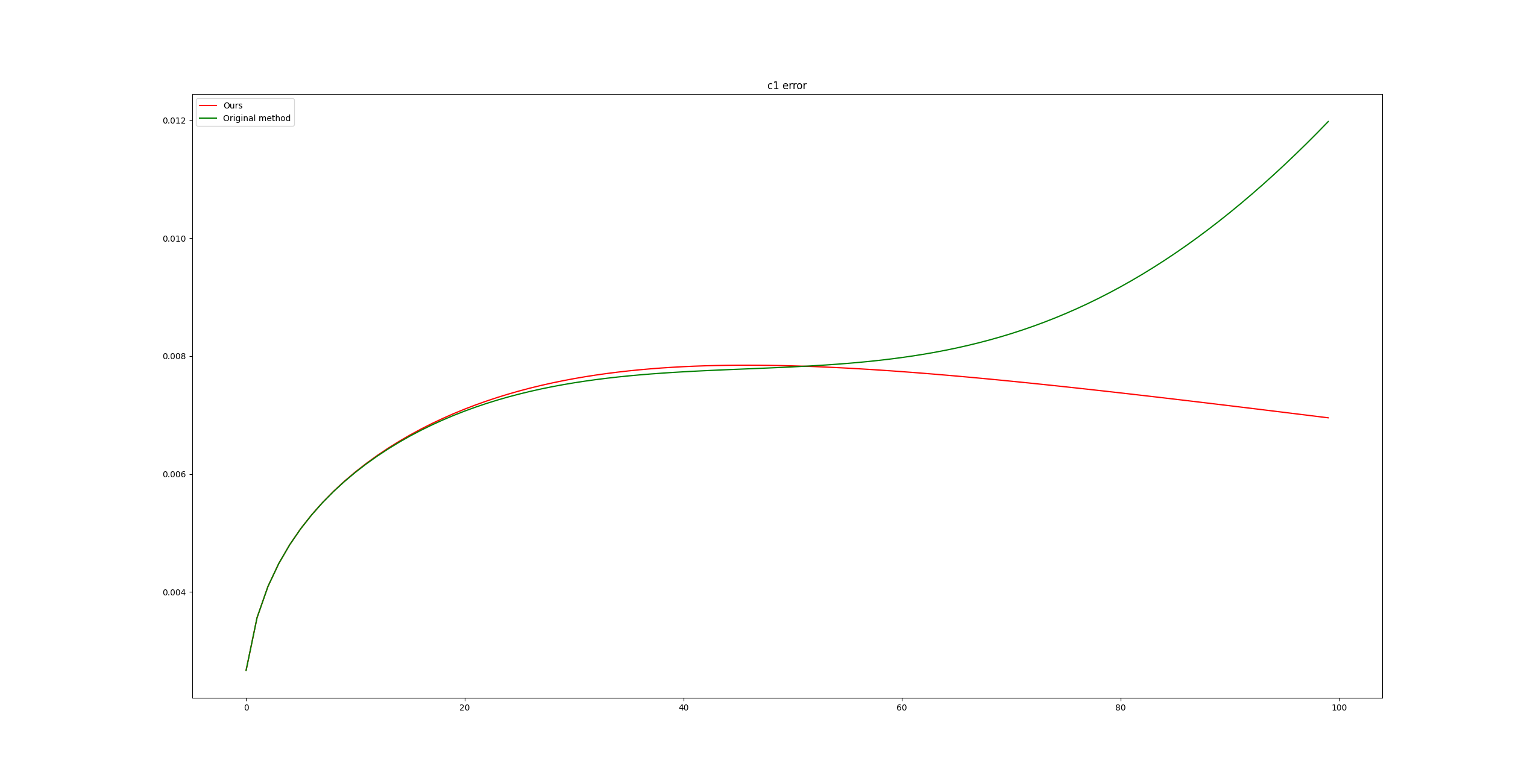}}
	\subcaptionbox{}{\includegraphics[scale=0.1]{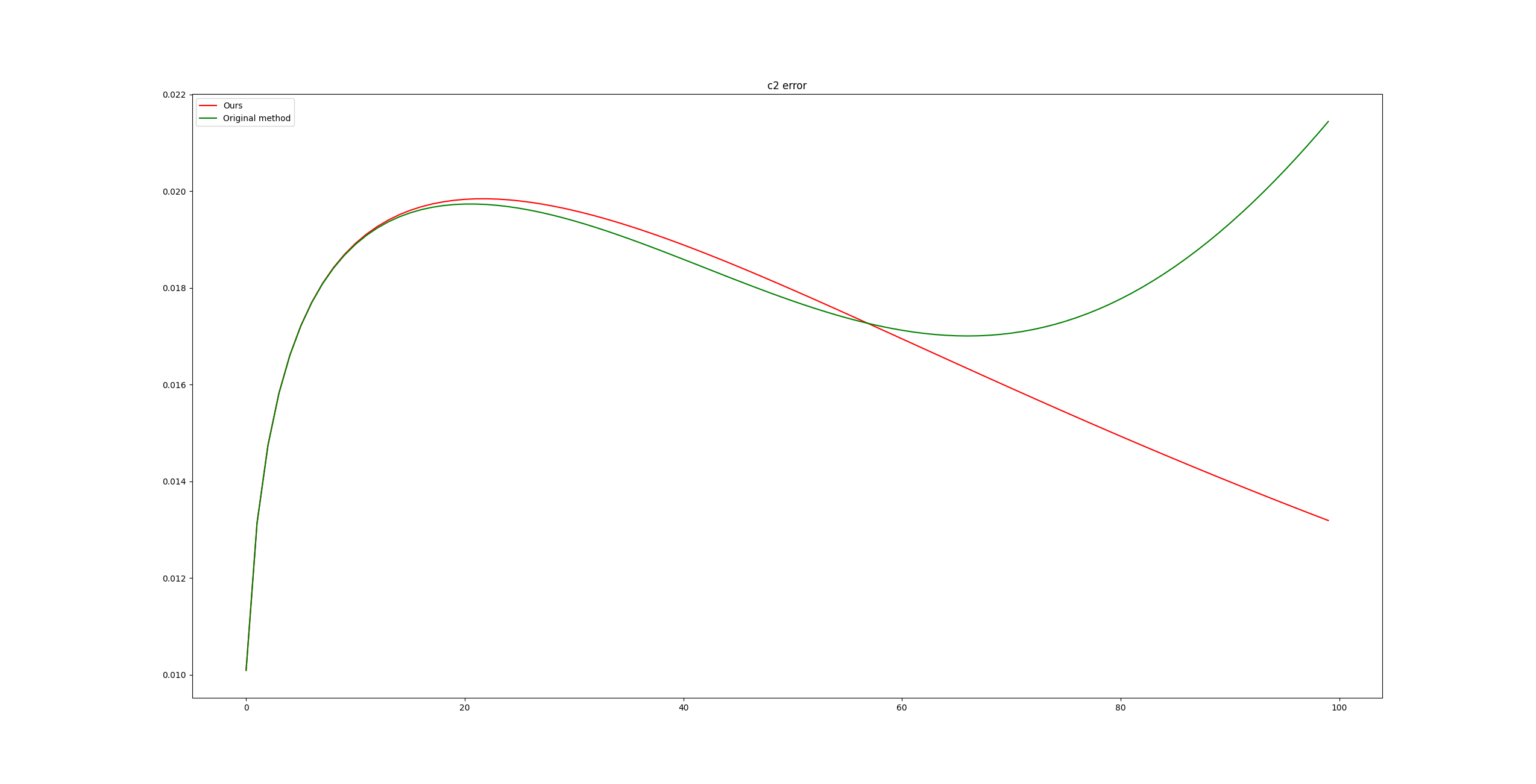}}
	\subcaptionbox{}{\includegraphics[scale=0.1]{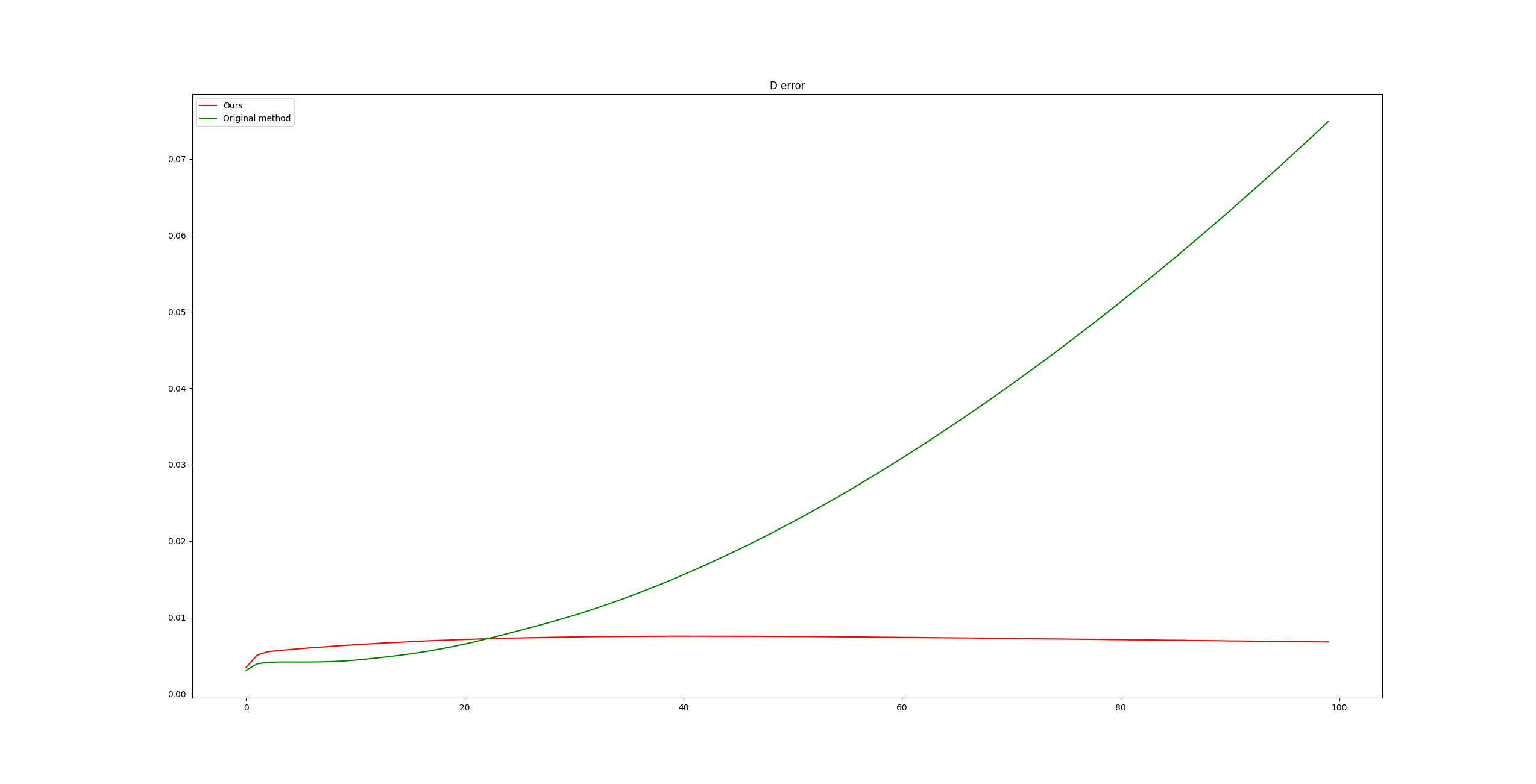}}
\centering
	\subcaptionbox{}{\includegraphics[scale=0.1]{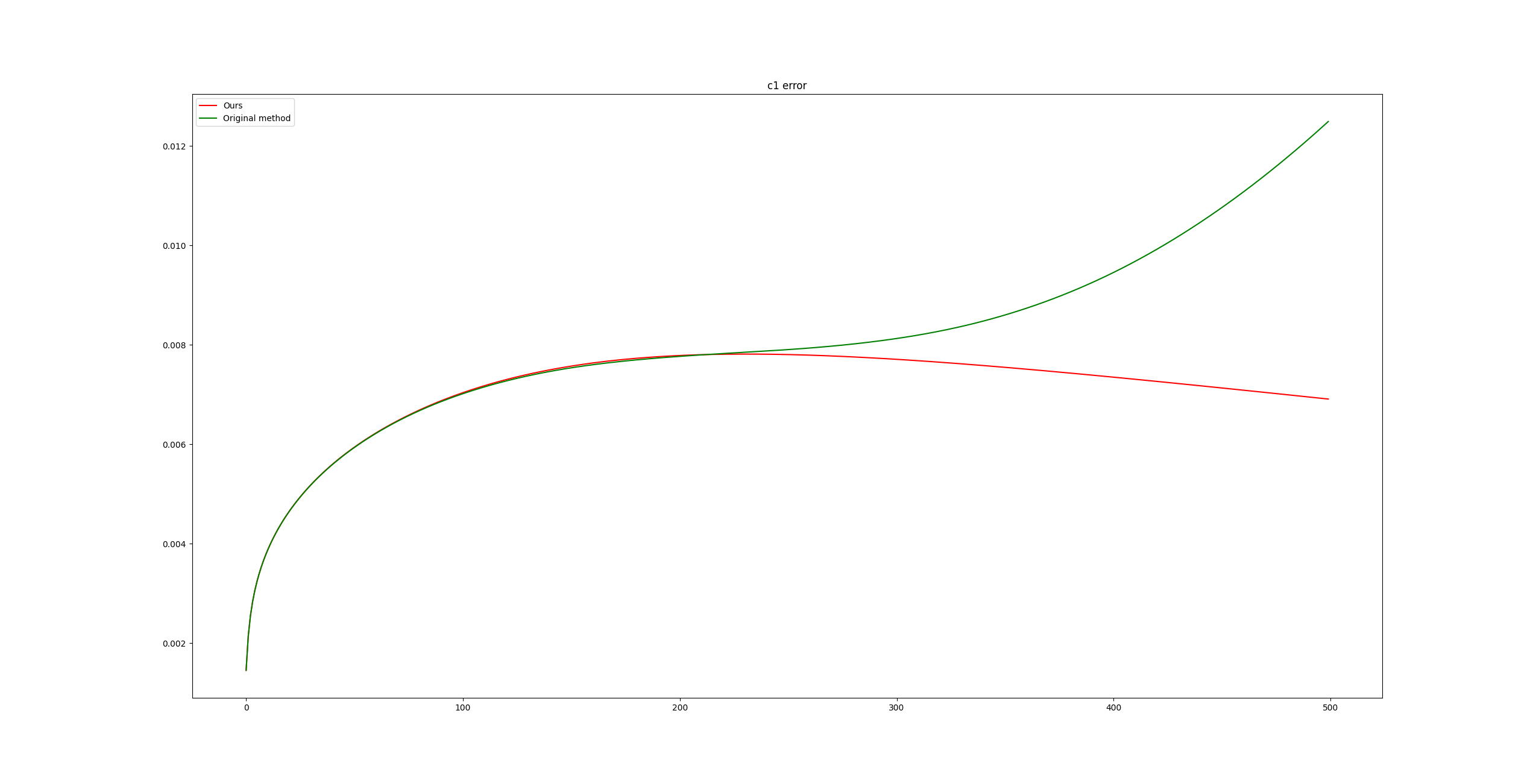}}
	\subcaptionbox{}{\includegraphics[scale=0.1]{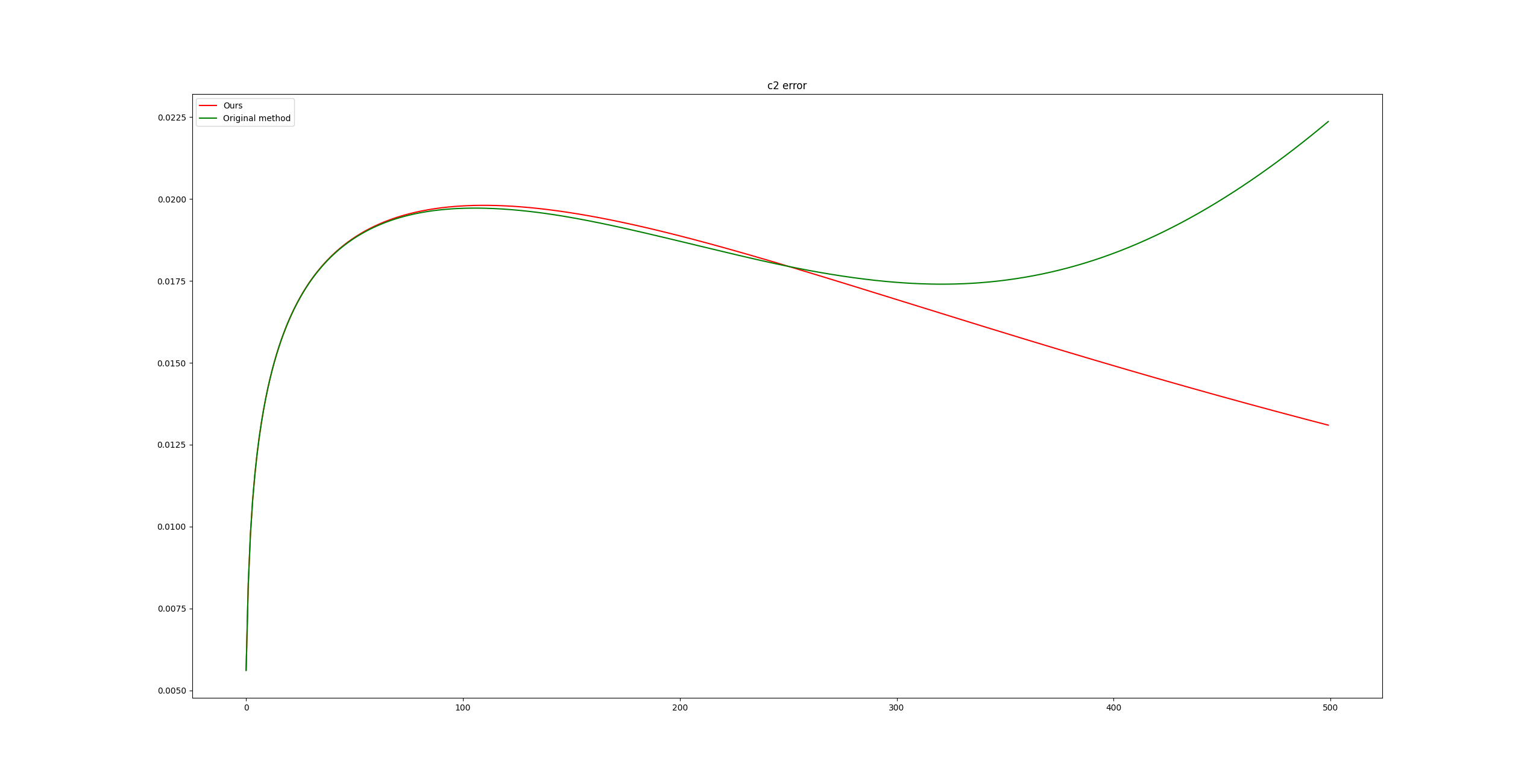}}
	\subcaptionbox{}{\includegraphics[scale=0.1]{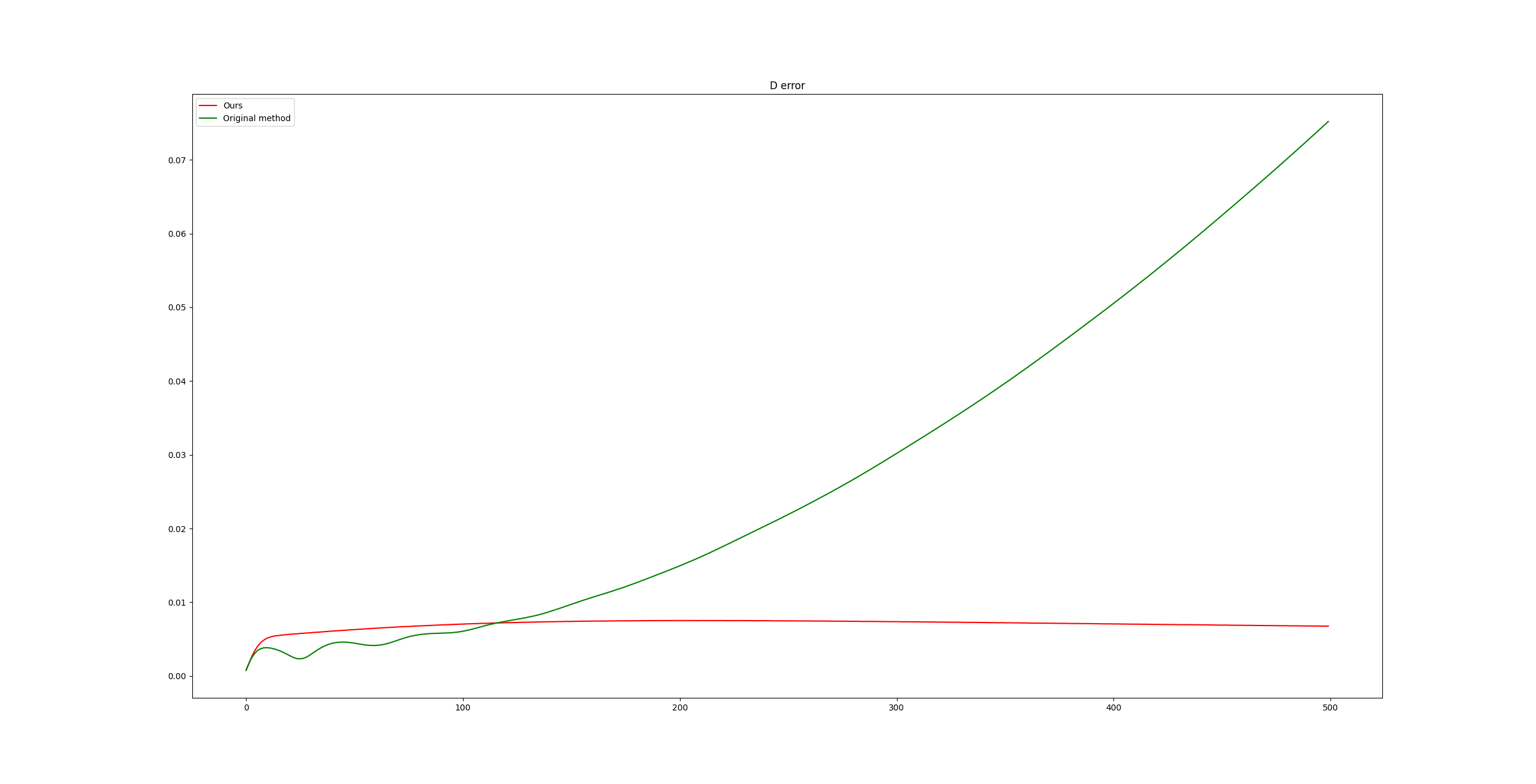}}
\caption{The error of $c^1$, $c^2$, and $\boldsymbol{D}$ at the different time steps in different experiment settings. Red lines are ours, and green lines are that of the original method \cite{ANP000_1, ANP000_2}. (a)--(c) $\Delta x=\Delta y=0.02$, $\Delta t=0.005$; (d)--(f) $\Delta x=\Delta y=0.02$, $\Delta t=0.001$; (g)--(i) $\Delta x=\Delta y=0.04$, $\Delta t=0.005$; (j)--(l) $\Delta x=\Delta y=0.04$, $\Delta t=0.001$}
\label{fig:error}
\end{figure*}

{We conduct several experiments with different resolutions, and record the error at each time step. In all experiments in this section and Section \ref{sec:2d_NE}, the acceleration method proposed in Section \ref{sec:Efficiency} is imposed on our method, while the original local curl-free algorithm is used in the original method as a baseline. The error of $c^1$, $c^2$ and $\boldsymbol{D}$ are shown in Figure \ref{fig:error}.}

{From the figure, it is clear that although our method shows comparable accuracy at the beginning, as time evolves, the error of our method still keeps constant or declines, while the that of the original method inclines. This phenomena can be observed in both $c^l$'s and $\boldsymbol{D}$, and is more significant in $\boldsymbol{D}$.}

\subsection{Electrodynamics in 2-dimensional space}\label{sec:2d_NE}
{The spatial domain is $\Omega=[-1, 1]\times[-1, 1]$, and the terminal time $T=0.5$. The cell size $\Delta x=\Delta y=0.04$, and $\Delta t=0.0005$.} We assume that the solute is a binary ionic compound with the valence of the cations and anions being $+1$ and $-1$, respectively. The compound is assumed to be fully ionised in the solvent, and the resultant ions can freely move in the solvent, under the effects of both themselves and fixed charges in the solvent. {The initial concentrations of both kinds of ions are 1 throughout the entire domain. $\epsilon$ is assumed to be 1. The fixed charges are described by the formula:}

\begin{equation}
{\rho^f=}
\left\{
\begin{array}{c}
{1, x\in\{(x-0.5)^2+y^2\leq 0.09\}}\\
{-1, x\in\{(x+0.5)^2+y^2\leq 0.09\}}\\
{0, \text{otherwise}}
\end{array}
\right.
\end{equation}

{The loss values of the neural network over time is shown in Figure \ref{fig:loss_2d}. A good convergence of the neural network is indicated by the stable downward trend of the loss values. In Figure \ref{fig:k_2d}, the number of iterations of the local curl-free relaxation algorithm \cite{ANP000_1, ANP000_2} is presented. The algorithm stops when the change of the objective function (Equation \ref{Local_curl_free_relaxation}) value is less than $10^{-5}$. Although the number of iterations is large at the very beginning, it drops sharply as the training progresses. It {remains 1 since the initial time step} (at least one iteration must be executed per time step to calculate the change of the objective function){, meaning that the neural network fit the appropriate value well, so that no further correction is required.}

\begin{figure}
  \centering
  \includegraphics[scale=0.4]{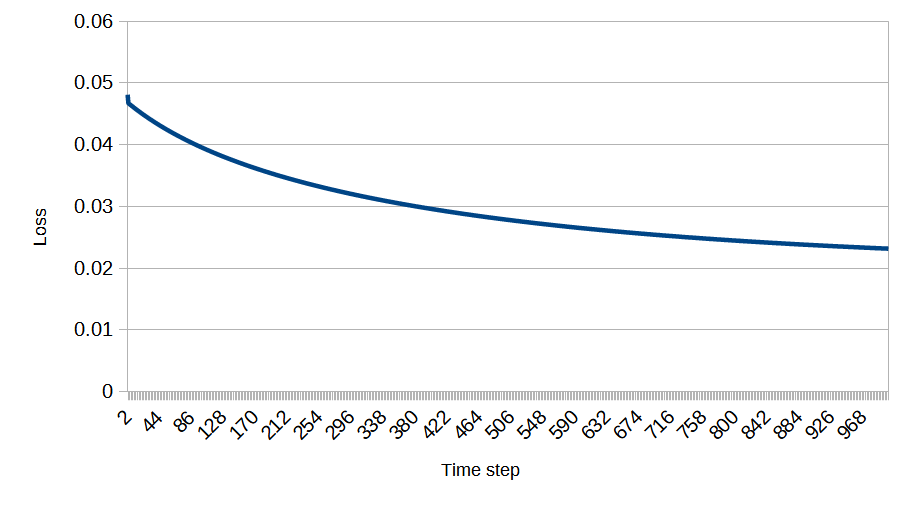}
  \caption{Neural network loss values vs. time steps.}
  \label{fig:loss_2d}
\end{figure}

\begin{figure}
  \centering
  \includegraphics[scale=0.4]{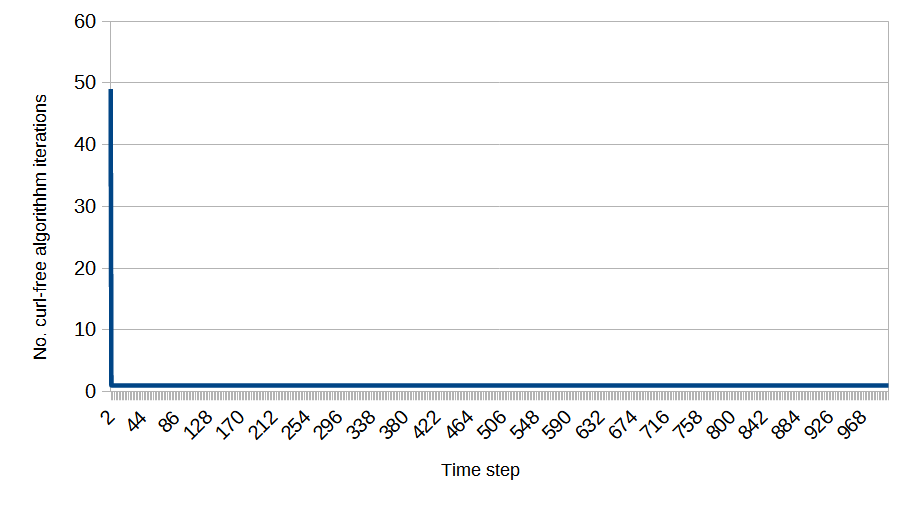}
  \caption{Number of iterations of local curl-free relaxation algorithm vs. time steps.}
  \label{fig:k_2d}
\end{figure}

The conservations of total mass of ions can be observed in Figure \ref{fig:c_total}. From Figure \ref{fig:c_min}, it can be seen that the positivity of ion concentrations {is} preserved. {Figure \ref{fig:energy} illustrates the evolution of the free energy of the system over time. Figure \ref{fig:c_D_norm} presents the snapshots of the numerical results at different time steps. A converging trend towards opposite charges can be observed from the ion concentration maps. At the borders of the fixed charges there is large electric displacement, and it fades as free ions move.}

\begin{figure*}
\centering
	\subcaptionbox{}{\includegraphics[scale=0.4]{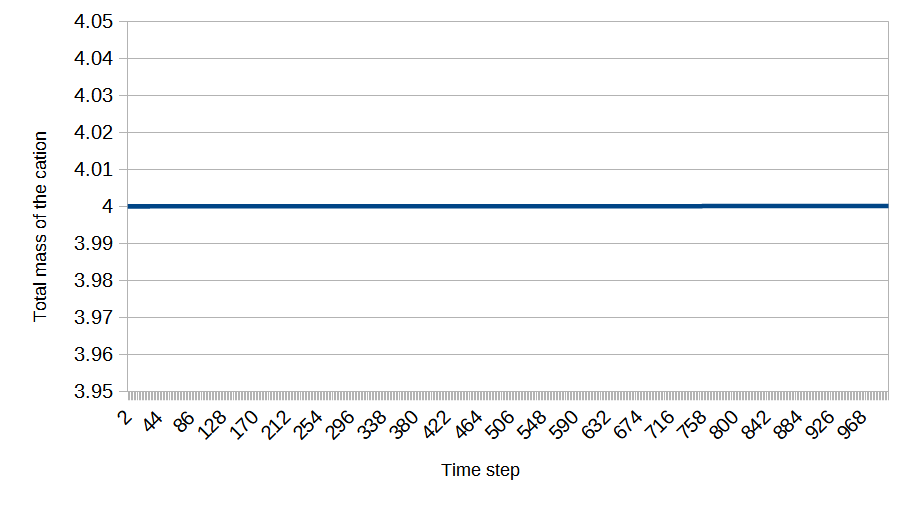}}
	\hfill
	\subcaptionbox{}{\includegraphics[scale=0.4]{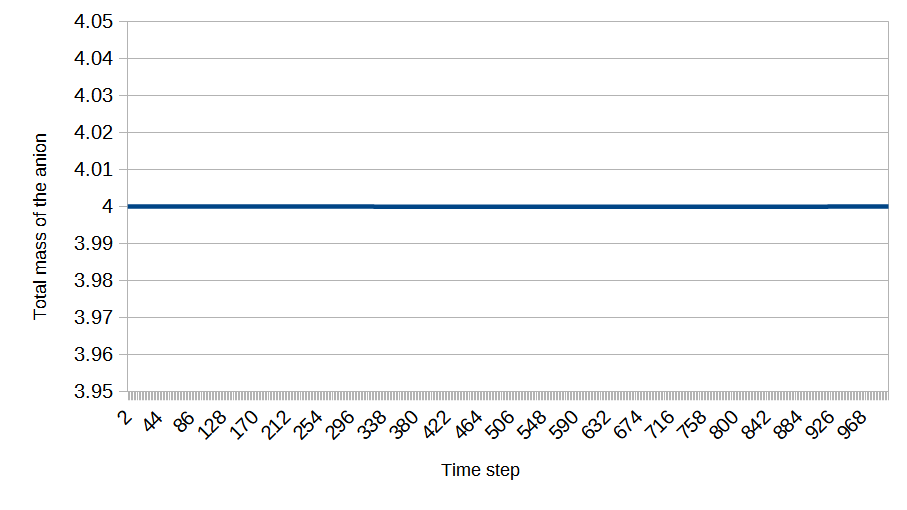}}
\caption{Total mass of the ions vs. time steps. (a) cation; (b) anion.}
\label{fig:c_total}
\end{figure*}

\begin{figure*}
\centering
	\subcaptionbox{}{\includegraphics[scale=0.4]{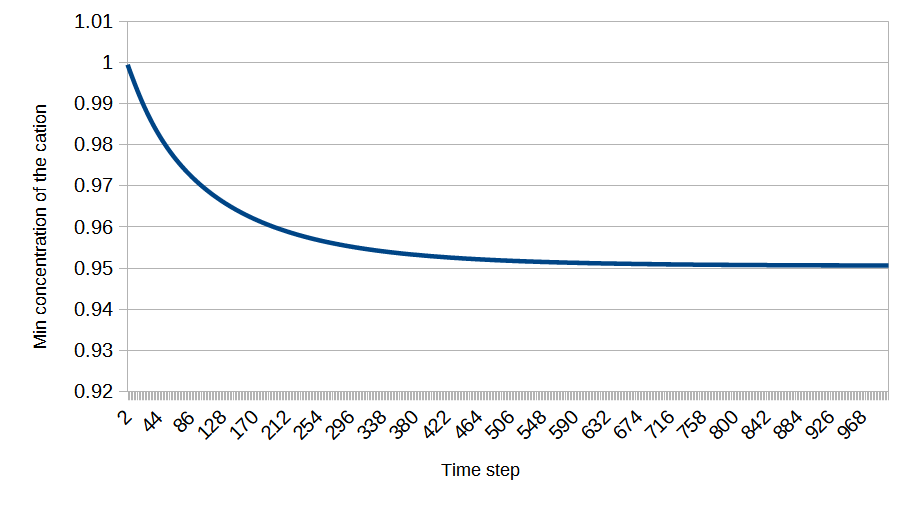}}
	\hfill
	\subcaptionbox{}{\includegraphics[scale=0.4]{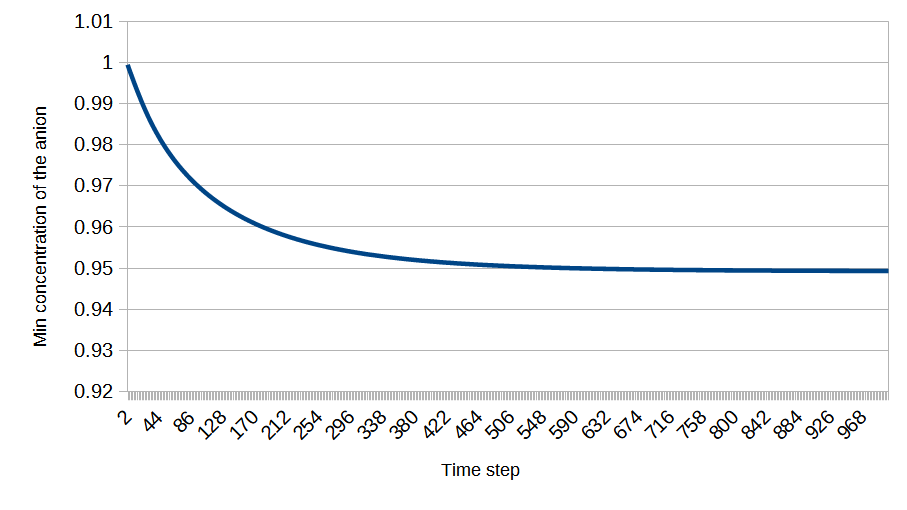}}
\caption{Minimum mass of cation vs. time steps. (a) cation; (b) anion.}
\label{fig:c_min}
\end{figure*}

\begin{figure}
  \centering
  \includegraphics[scale=0.4]{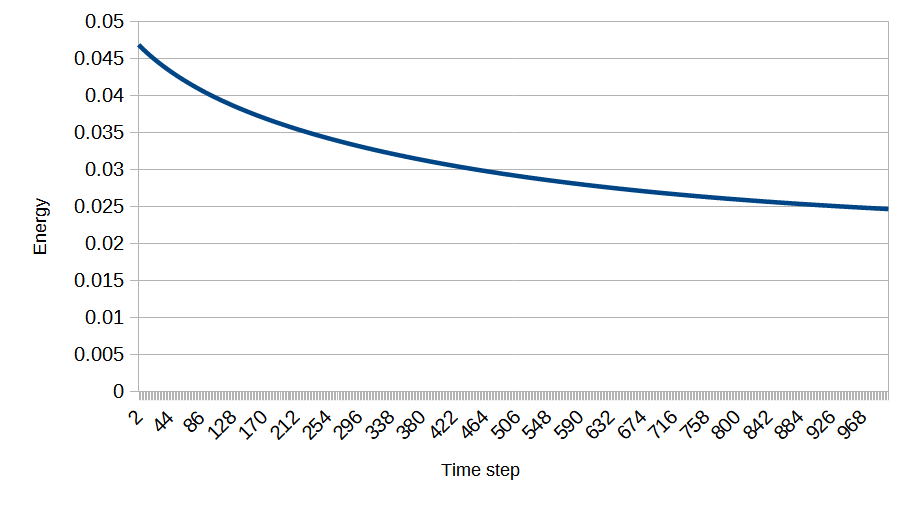}
  \caption{Free energy vs. time steps.}
  \label{fig:energy}
\end{figure}

\begin{figure*}
\centering
	\subcaptionbox{}{\includegraphics[scale=0.3]{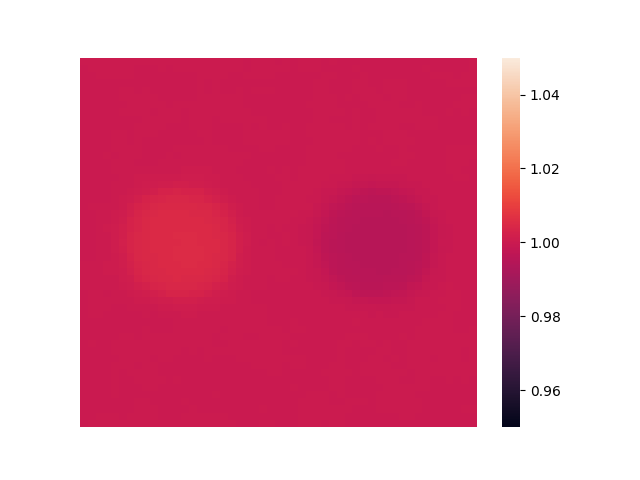}}
	\subcaptionbox{}{\includegraphics[scale=0.3]{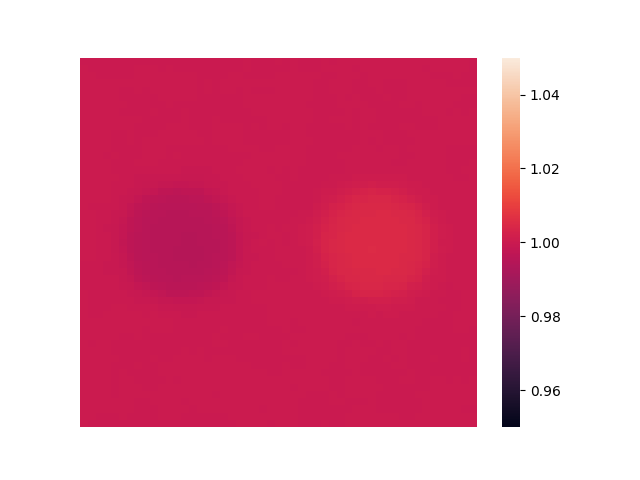}}
	\subcaptionbox{}{\includegraphics[scale=0.3]{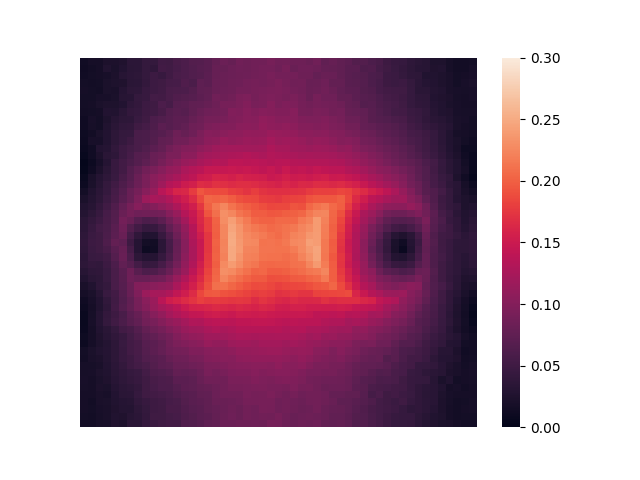}}
\centering
	\subcaptionbox{}{\includegraphics[scale=0.3]{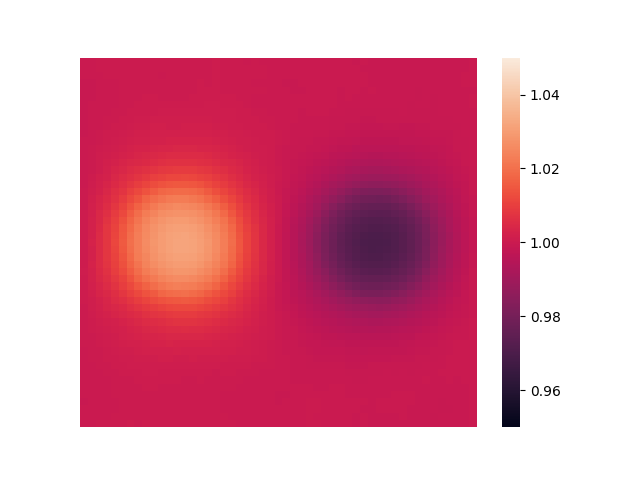}}
	\subcaptionbox{}{\includegraphics[scale=0.3]{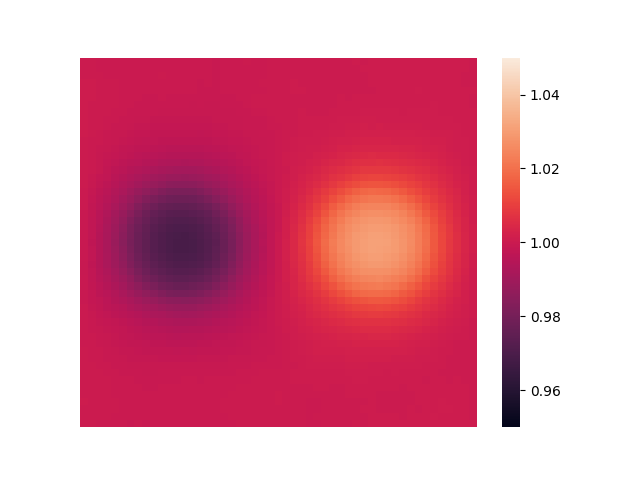}}
	\subcaptionbox{}{\includegraphics[scale=0.3]{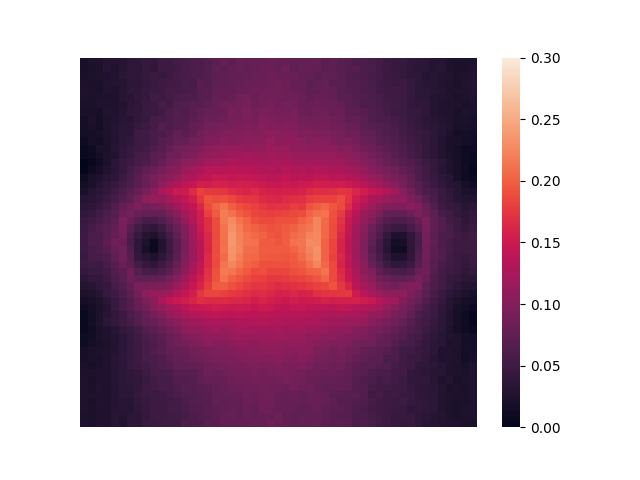}}
\centering
	\subcaptionbox{}{\includegraphics[scale=0.3]{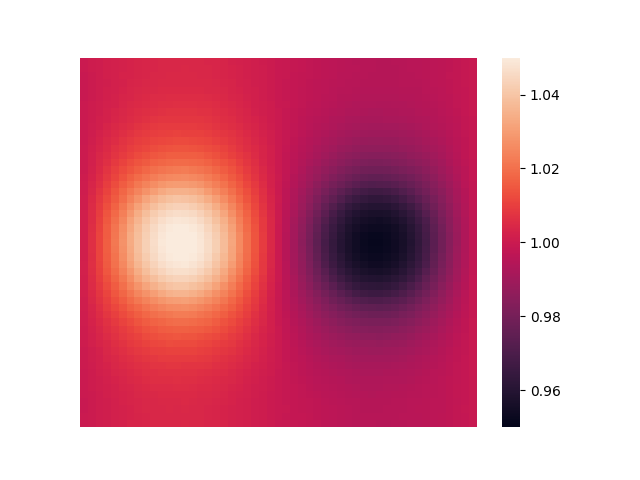}}
	\subcaptionbox{}{\includegraphics[scale=0.3]{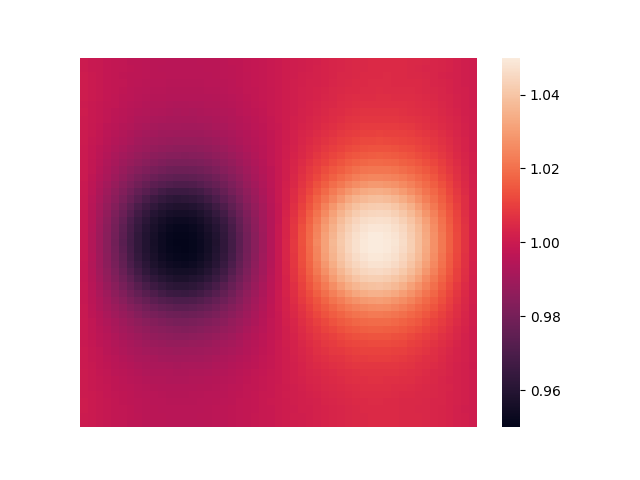}}
	\subcaptionbox{}{\includegraphics[scale=0.3]{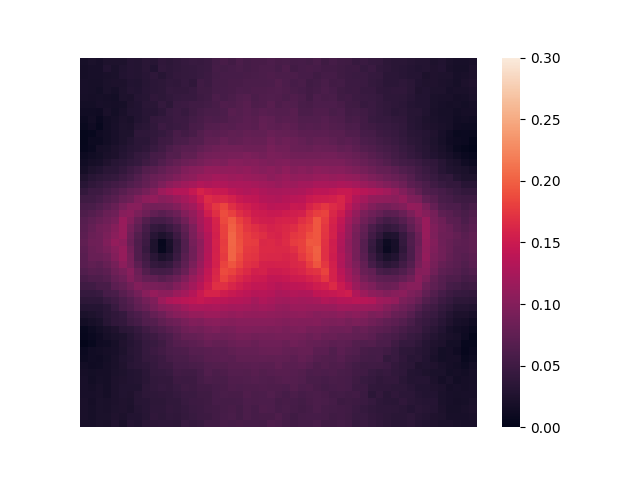}}
\centering
	\subcaptionbox{}{\includegraphics[scale=0.3]{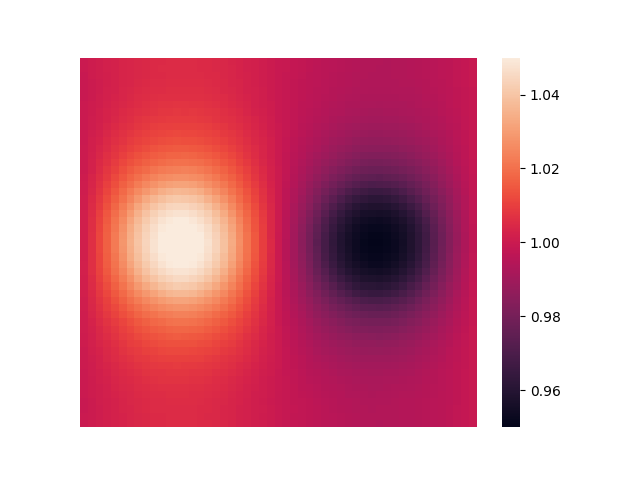}}
	\subcaptionbox{}{\includegraphics[scale=0.3]{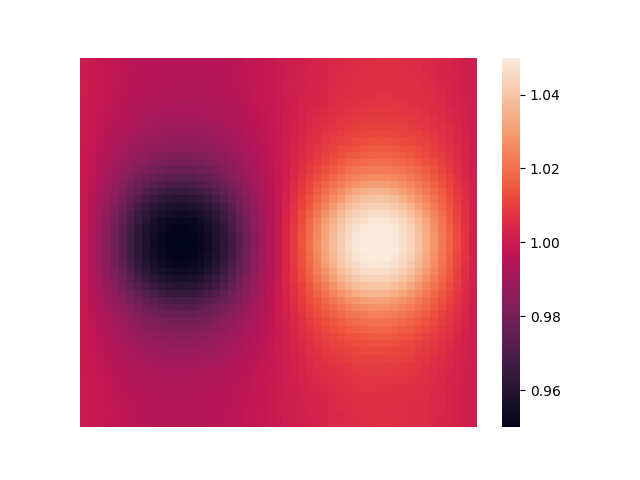}}
	\subcaptionbox{}{\includegraphics[scale=0.3]{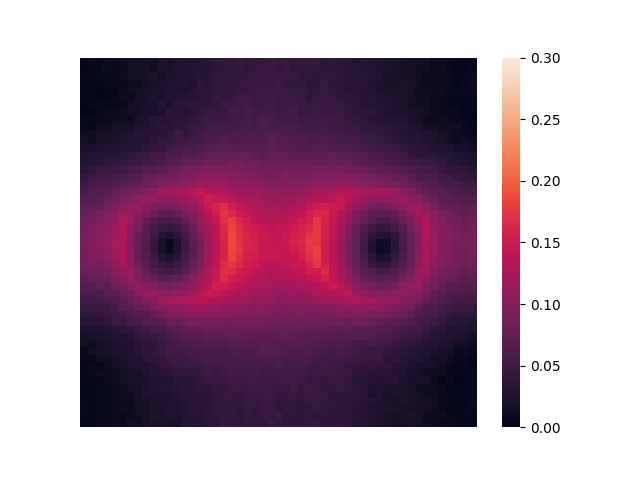}}
\caption{$c^1$, $c^2$, and the $l_2$-norm of $\boldsymbol{D}$ at the different time steps. (a)--(c) Results at time step 10; (d)--(f) results at time step 100; (g)--(i) results at time step 500; (j)--(l) results at time step 1000}
\label{fig:c_D_norm}
\end{figure*}

\section{Conclusions}

{This study proposes} a hybrid numerical method combining{ }conventional methods and deep learning approaches {to solve} the {M}ANP equations. {A} deep neural network {is employed} to improve {the} conventional numerical scheme. The resulting hybrid method combines the strengths of both {conventional and deep learning methods}. First, it{ }preserves the conservation properties of the {M}ANP equations with particular boundary condition{s}, including the conservation of {the} total mass of ions {and} the positivity of the {ion concentrations.} Second, {owing} to the flexibility and universal applicability of deep learning, the proposed method exhibits good portability to different problem settings. {Numerical} experiments on both 1-dimensional and 2-dimensional cases are conducted to validate the method. They {demonstrates} the generalisation capability of the proposed method{ }as well as the preservation of conservation properties{ }.

\bibliographystyle{unsrt}  

\bibliographystyle{siam}
\bibliography{references.bib}

\end{document}